\documentclass{imanum}

\usepackage{latexsym,amssymb,amscd,amsmath}
\usepackage{amsfonts,enumerate,epsfig}
\usepackage{latexsym,amssymb,amscd,amsmath,graphics,color,makeidx}

\gridframe{N}

\jno{dri017}

\usepackage{modamsthm,modnatbib}

\usepackage{latexsym,amssymb}
\usepackage{lastpage}
\usepackage{graphicx,amsfonts}
\usepackage{times,mathptmx,bm,amsmath}
\usepackage{dcolumn}

 \newtheoremstyle{theorem}{6pt}{6pt}{\rm}{}{\sffamily}{ }{ }{}
 \theoremstyle{theorem}

  \newtheoremstyle{thm}{6pt}{6pt}{\rm}{}{\sffamily}{ }{ }{}
 \theoremstyle{thm}

 \newtheoremstyle{lemma}{6pt}{6pt}{\rm}{}{\sffamily}{ }{ }{}
 \theoremstyle{lemma}
 
 \newtheoremstyle{lem}{6pt}{6pt}{\rm}{}{\sffamily}{ }{ }{}
 \theoremstyle{lem}

\newtheoremstyle{case}{6pt}{6pt}{\rm}{}{}{. }{ }{}
 \theoremstyle{case}

 \newtheoremstyle{statement}{6pt}{6pt}{\rm}{}{\sffamily}{ }{ }{}
\theoremstyle{statement}

 \newtheoremstyle{corollary}{6pt}{6pt}{\rm}{}{\sffamily}{ }{ }{}
 \theoremstyle{corollary}

  \newtheoremstyle{defi}{6pt}{6pt}{\rm}{}{\sffamily}{ }{ }{}
 \theoremstyle{defi}

  \newtheoremstyle{cor}{6pt}{6pt}{\rm}{}{\sffamily}{ }{ }{}
 \theoremstyle{cor}

\newtheoremstyle{example}{6pt}{6pt}{\rm}{}{\sffamily}{ }{ }{}
\theoremstyle{example}

\newtheoremstyle{remark}{6pt}{6pt}{\rm}{}{\sffamily}{ }{ }{}
\theoremstyle{remark}

\newtheoremstyle{approximation}{6pt}{6pt}{\rm}{}{\sffamily}{ }{ }{}
\theoremstyle{approximation}

\newtheoremstyle{scheme}{6pt}{6pt}{\rm}{}{\sffamily}{ }{ }{}
\theoremstyle{scheme}

\newtheoremstyle{Algorithm}{6pt}{6pt}{\rm}{}{\sffamily}{ }{ }{}
\theoremstyle{Algorithm}

 \newtheoremstyle{Remark}{6pt}{6pt}{\rm}{}{\sffamily}{ }{ }{}
 \theoremstyle{Remark}

\newtheoremstyle{Lemma}{6pt}{6pt}{\rm}{}{\sffamily}{ }{ }{}
\theoremstyle{Lemma}

\newtheoremstyle{Assumption}{6pt}{6pt}{\rm}{}{\sffamily}{ }{ }{}
\theoremstyle{Assumption}

\newtheoremstyle{Proposition}{6pt}{6pt}{\rm}{}{\sffamily}{ }{ }{}
\theoremstyle{Proposition}

\newtheoremstyle{prop}{6pt}{6pt}{\rm}{}{\sffamily}{ }{ }{}
\theoremstyle{prop}
\newtheorem{prop}{\sc Proposition}[section]

\newtheoremstyle{rem}{6pt}{6pt}{\rm}{}{\sffamily}{ }{ }{}
 \theoremstyle{rem}
 \newtheorem{rem}{\sc Remark}[section]

\newtheoremstyle{hypo}{6pt}{6pt}{\rm}{}{\sffamily}{ }{ }{}
 \theoremstyle{hypo}

  \newtheoremstyle{Step}{6pt}{6pt}{\rm}{}{}{ }{ }{}
 \theoremstyle{Step}

 \newtheoremstyle{lema}{6pt}{6pt}{\rm}{}{\sffamily}{ }{ }{}
 \theoremstyle{lema}

\renewcommand{\theequation}{\thesection.\arabic{equation}}
\numberwithin{equation}{section}

\def \qed {\hfill \rule{1ex}{1ex}}
\def \U {{\bf u}}
\def \x {{\bf x}}

\def \H {{\bf H}}

\def \V {{\bf v}}
\def \bt {{\widetilde{\hbox{\rm b}}}}
\def \W {{\bf w}}
\def \N {{\mathbf{n}}}

\def \Ut {{\tilde {\mathbf{u}}}}

\def \F {{\bf f }}
\def \P {{\bf P}}

\def \div {{\hbox{div}}}

\def \L {{\bf L}}

\def \b {{\hbox{\rm b}}}

\let \eps=\varepsilon

\newcommand{\Eta}[2]{\hbox{\boldmath $ \eta^{\hbox{\unboldmath $\scriptstyle#1$}}
_{\hbox{\unboldmath $\scriptstyle #2$}}$ \unboldmath \hspace{-3mm}
}}

\newcommand{\Eps}[2]{\hbox{\boldmath $ \epsilon^{\hbox{\unboldmath $\scriptstyle#1$}}
_{\hbox{\unboldmath $\scriptstyle #2$}}$ \unboldmath
\hspace{-1.6mm}}}

\newtheorem{theor}{Theorem}[section]

\begin{document}
\title{Study of a finite volume- finite element  scheme for a nuclear transport model}
\author{{\sc Catherine Choquet }\\[2pt]
LATP, Universit\'e Aix-Marseille 3, 13013 Marseille Cedex 20, France \\[6pt]
 {\sc S\'ebastien Zimmermann}\\[2pt]
 Laboratoire de math\'ematiques, \'Ecole Centrale de Lyon, 36 avenue Guy de Collongue, 69134 Ecully, France. \\[6pt]
{\rm [Received on  24 July 2007]}} \pagestyle{headings} \markboth{C.
CHOQUET AND S. ZIMMERMANN}{\rm Study of a finite volume-finite
element  scheme for a nuclear transport model} \maketitle


\begin{abstract}
{ We consider  a problem of nuclear waste contamination.
It  takes into account the thermal effects. The temperature and the
contaminant's concentration  fulfill convection-diffusion-reaction
equations. The velocity and the pressure  in the flow satisfy  the
Darcy equation, with a viscosity depending on both concentration and
temperature. The equations are nonlinear and strongly coupled. Using
both finite volume and nonconforming finite element methods, we
introduce a scheme adapted to this problem. We prove the stability
and convergence of this scheme and give some error estimates. }
{porous media, miscible flow, nonconforming finite element, finite
volume}.
\end{abstract}

\section{Introduction}

A part of the high-level nuclear waste is now stored in
environmentally safe locations. One has to consider the eventuality
of a leakage through the engineering and geological barriers. It may
cause the contamination of underground water sources far away from
the original repository's location. In the present paper, we
consider such  a problem of nuclear waste contamination in the
basement.
We take into account the thermal effects. The evolution in time of
the temperature and of the contaminants concentration  is then
governed by convection-diffusion-reaction equations. The velocity
and the pressure  in the flow satisfy  the Darcy equation, with a
viscosity depending on both concentrations and temperature in a
nonlinear way. The velocity satisfies an incompressibility
constraint.
We introduce a scheme adapted to  this problem. We use both finite
volume and nonconforming finite element methods. It ensures that a
maximum principle holds and that the associated linear systems have
good-conditioned matrices. We prove the stability and convergence of
the scheme and give some error estimates.

 Let us briefly point out some previous works. A complete model
coupling concentrations and pressure equations is very seldom
studied, since the system is strongly coupled. Instead, each
equation is considered separately. In the general context of
convection-diffusion-reaction equations, numerous schemes are
available (see  \cite{john}
or  \cite{morton}
and the references therein). Finite difference schemes are sometimes
used for the convective term (in \cite{verw} for instance). But they
are not adapted for the complex geometry of a reservoir. More
recently, finite volume methods were developed and analysed. Let us
just cite the book \cite{afif} or \cite{eym}
and the references therein. Finite elements (for the diffusive term)
and finite volumes (for the convective term) are coupled for
instance in \cite{angot,chav2}
. In convection dominant problems, the equations are of degenerate
parabolic type. This setting is considered in \cite{chen1,mich}
. The reaction terms are specifically studied through operator
splitting methods in \cite{geis}
. Now, in the specific context of porous media flow, we mention
\cite{lin} who consider only the evolution of the  pressure.
  In \cite{doug,ewing2} a more complete set of equations is
  used,
and  a mixed finite element approximation is developed. 
We stress that in all theses works, as in most, the
thermic effects are not taken into account.

 The present paper is organized as follows. Section \ref{sec:model}
is devoted to the derivation of the model. In section
\ref{sec:dtools}, we introduce the discrete tools used in this
paper. It allows us to define the numerical scheme  of section
\ref{sec:schema}. The analysis of the scheme uses the properties of
section \ref{sec:propop}. We then prove the stability and
convergence of the scheme,  in  sections \ref{sec:stab} and
\ref{sec:conv} respectively. We conclude with some error estimates
in section \ref{sec:esterr}.

\section{Model of contamination}
\label{sec:model}

\newcounter{matmod}
\renewcommand{\theequation}{\thesection.\thematmod}

The thickness of the medium is significantly smaller than its length
and width. Hence it is reasonable to average the medium properties
vertically and to describe the far-field repository by a polyhedral
domain $\Omega$ of $\mathbb{R}^2$  with a smooth boundary $\partial
\Omega$. It is characterized by a porosity $\phi$ and a permeability
tensor $K$.
The time interval of interest is $[0,T]$. We denote by $p$ the
pressure, by $(a_i)_{i=1}^{N_r}$ the concentrations of the $N_r$
radionuclides involved in the flow and by $\theta$ the temperature.
The Darcy velocity is represented by $\U$. We assume a miscible and
incompressible displacement. Due to the mass and energy conservation
principles, the flow is governed by the following system satisfied
in $\Omega \times [0,T]$, with $i =1,..,N_r-1$ (see
\cite{choquet}). 
 \stepcounter{matmod}
\begin{eqnarray}
&& \phi  \, R_i \, \partial_t a_i + \div(a_i  \, \U ) -  \div (\phi
\, D_c \, \nabla a_i) =s_i  - s \, a_i -\lambda_i\,  R_i\,  \phi\,
a_i+\sum_{j=1,j \ne i}^{N_r-1} k_{ij}\,  \lambda_j\,  R_j\,  \phi\,
a_j \, ,
\label{eq:1.1}\\
\stepcounter{matmod} && \phi \, \mathcal{C}_p \, \partial_t \theta +
\div (\theta \, \U) -  \div ( \phi  \, D_\theta \, \nabla \theta) =
- s_\theta  - s \, (\theta-\theta_*) \, ,
\label{eq:1.2} \\
\stepcounter{matmod} && \div \, \U = s, \quad \U
+\frac{K}{\mu\big((a_j)_{j=1}^{N_r-1}, \theta\big)}\,  \nabla p =
\F. \label{eq:1.3}
\end{eqnarray}
In (\ref{eq:1.1}) the retardation factors $R_i
>0$ are due to the sorption mechanism. The real $\lambda_i^{-1} >0$
denotes the half life time of radionuclide $i$. The term $-\lambda_i
\, R_i \, \phi \,  a_i$ describes the radioactive decay of the
i--{\it th} specy. Meanwhile, the quantity $\sum_{j \ne i} k_{ij} \,
\lambda_j \, R_j \, \phi \, a_j$ is created by radioactive
filiation. The molecular diffusion effects are given by the
coefficient $D_c>0$. The contamination is represented by the  source
term $s_i$ and
 $s=\sum_{i=1}^{N_r} s_i$.
 In (\ref{eq:1.2}) the coefficient $\mathcal{C}_p>0$ is the relative specific heat of the
porous medium. The thermic diffusion coefficient is denoted by
$D_\theta>0$.
The real $\theta_*>0$ is a reference temperature. The constitutive
relation (\ref{eq:1.3}) is the Darcy law and $\F$ is a density of
body forces. For a large range of temperatures $\mu$ has the form
$$\mu(a,\theta)= \mu_R(a) \exp \left(
    \frac{1}{\theta} - \frac{1}{\theta_*}  \right)   $$
where $\mu_R$ is a nonlinear function. For instance, in the Koval
model for a two-species mixture \cite{koval}, we have
$$\mu_R(a)=\mu(0)(1+(M^{1/4}-1)a_1)^{-4}$$
 where $M=\mu(0)/\mu(1)$ is the mobility ratio.

\noindent We notice that the equations (\ref{eq:1.1})-(\ref{eq:1.3})
are strongly coupled. Moreover,
  every concentration equation (\ref{eq:1.1}) involves a  different time scale. Therefore, it is difficult to build a numerical scheme that captures all the physical phenomena. We have to transform these equations.
We first assume that only serial or parallel first-order reactions
occur, so that $k_{ij}=y_i$ with $y_1=0$.
Next, following \cite{bate}, we assume that no two isotopes have
identical decay rates and we set
\stepcounter{matmod}
\begin{equation}
c_1=a_1, \qquad c_i = a_i + \sum_{j=1}^{i-1} \left(
\prod_{l=j}^{i-1} \frac{y_{l+1} \lambda_l}{\lambda_l - \lambda_i}
\right)  a_j \quad \mathrm{for}\  i=2,..,N_r . \label{eq:1.9}
\end{equation}
Lastly, without losing any mathematical difficulty (see remark
\ref{rem:rem2}  below), we set $R_i=1$ for $i=1,..,N_r-1$ and
$\phi=1$, ${\mathcal{C}}_p=1$.
We also set  $
s_{c_i} = s_i + \sum_{j=1}^{i-1} \bigl(\prod_{l=j}^{i-1}
\frac{y_{l+1}\lambda_l}{\lambda_l - \lambda_i} \bigr) \, s_j
 $ and $\kappa(c,\theta)=K/\mu\big((a_i)_{i=1}^{N_r-1},\theta\big)$.
The contamination problem is now modelized by the following
parabolic-elliptic system \stepcounter{matmod}
\begin{eqnarray}
&&  \partial_t c_i +  \div(c_i \, \U ) -  D_c \, \Delta c_i =
s_{c_i} -s \, c_i-\lambda_i\,  c_i \, ,
\label{eq:1.10} \\
\stepcounter{matmod} && \partial_t \theta + \div (\theta \, \U) -
D_\theta \, \Delta \theta = - s_\theta -s \, (\theta-\theta_*) \, ,
\label{eq:1.11} \\
\stepcounter{matmod} && \div \, \U = s, \quad \U +\kappa(c,\theta)
\,  \nabla p = \F \, , \label{eq:1.12}
\end{eqnarray}
with  $i=1,..,N_r-1$. These equations are completed  with the
boundary and initial conditions \stepcounter{matmod}
\begin{eqnarray}
 \label{eq:1.13}
  && \nabla c_i \cdot \nu =0 \, , \hspace{1cm}
 \nabla \theta \cdot \nu =0 \, , \hspace{1cm} \U \cdot \N =0 \, , \\
\stepcounter{matmod}
&& \label{eq:1.14}
 c_i|_{t=0} = c^i_{0}, \hspace{1cm} \theta|_{t=0}=\theta_0.
\end{eqnarray}
The pressure $p$ is normalized by $\int_\Omega p \, d\x=0$.
Equations (\ref{eq:1.10}) are all similar. Thus, for the sake of
simplicity, we will assume that there is only one. We set $N_r=2$
and $c=c_1$, $c_0=c^1_{0}$, $s_c=s_{c_1}$, $\lambda=\lambda_1$. The
results of this paper readily extend to the general case.


\noindent We conclude with some notations and hypothesis. Let  $D$
be a bounded open set of $\mathbb{R}^k$ with $k \ge 1$. We denote by
${\mathcal{C}}^\infty_0(D)$ the set of functions that are continuous
on $D$ together with all their derivatives, and have  a compact
support in $D$. For $p\in\{2,\infty\}$, we use the Lebesgue spaces
$\big(L^p(D),\|.\|_{L^p(D)}\big)$ and
$\big(\L^p(D),\|.\|_{\L^p(D)}\big)$ with $\L^p=(L^p)^2$. We also use
the Sobolev spaces $W^{p,q}(D)$ for $p\in[1,\infty[$ and
$q\in[1,\infty[$.
 In the case $D=\Omega$ we use the following conventions. We drop the domain dependancy.
We denote by $|.|$ (resp. $\|.\|_\infty$) the norms associated to
$L^2=L^2(\Omega)$ and $\L^2=\L^2(\Omega)$ (resp. $L^\infty$ and
$\L^\infty$). We set $L^2_0= \{ v \in L^2  ; \   \int_\Omega v(\x)
\, d\x=0\}$. For $p\in[1,\infty[$ we define $(H^p,\|.\|_p)$ and
$(\H^p,\|.\|_p)$ with $H^p=W^{p,2}$  and $\H^p=(H^p)^2$. Now let
$(X,|.|)$ be  a Banach space. The  functions $g:[0,T] \to X$ such
that $t \to \|g(t)\|_X$ is continuous (resp. bounded and square
integrable) form the set ${\mathcal{C}}(0,T;X)$ (resp.
$L^\infty(0,T;X)$ and $L^2(0,T;X)$). The associated norm for the
space  $L^\infty(0,T;X)$ (resp. $L^2(0,T;X)$) is defined by
$\|g\|_{L^\infty(0,T;X)}=\sup_{t\in [0,T]} \|g(t)\|_X$ (resp.
 $\|g\|_{L^2(0,T;X)}=\left(\int_{0}^T \|f(t)\|^2_X \,
dt\right)^{1/2}$). Finally in all  computations we use  $C>0$ as a
generic constant. It depends only on the data
 of the problem.

\noindent We assume the following regularities for the data in
(\ref{eq:1.10})--(\ref{eq:1.12}) \stepcounter{matmod}
\begin{equation}
\label{eq:regdata}
  \kappa \in W^{1,\infty}((0,1)\times(0,\infty)) \,, \hspace{1cm} s ,\, s_c , \, s_\theta \in L^2 \, , \hspace{1cm}
  \F \in {\mathcal{C}}(0,T;\L^2).
\end{equation}
We also assume that $\kappa \ge \kappa_{inf}$ with $\kappa_{inf}>0$.
For the initial data, we assume that $c_0 \in H^1$, $\theta_0 \in
H^1$, and that we have a.e. in $\Omega$ \stepcounter{matmod}
\begin{equation}
\label{eq:borncondi}
  0 \le c_0(\x) \le 1 \, , \hspace{2cm} \theta_- \le \theta_0(\x) \le \theta_+
\end{equation}
with $\theta_->0$. Finally we assume that  we have a.e. in $\Omega$
\begin{equation}
\stepcounter{matmod} \label{eq:condpmax}
  2 \, s(\x) +\lambda \ge s_c(\x) \ge 0 \, , \hspace{.6cm}
  2 \, (\theta^{-} -\theta^*) \, s(\x)+s_\theta(\x) \le 0 \, ,
  \hspace{.6cm} 2 \, (\theta^{+} -\theta^*) \, s(\x)+s_\theta(\x) \ge 0.
\end{equation}
These conditions ensure a maximum principle (proposition
\ref{propmaxd} below).

\begin{rem}
\label{rem:rem1} We have assumed that first-order reactions occur,
so that the coefficients  $k_{ij}$ in (\ref{eq:1.1}) depend only on
$i$. If $k_{ij}$ depends on $i$ and $j$, one can still uncouple the
equations by iterating  the transformation (\ref{eq:1.9}), provided
that $k_{1j}=0$ for $j=2,..,N_r-1$. This assumption  means that the
first long-lasting isotope disappears and is not created anymore. It
is satisfied by many radionuclides.
\end{rem}

\begin{rem}
\label{rem:rem2} We have assumed that the retardation factors $R_j$
are identical. If it is not the case, the difficulty and the
approach remain the same. Indeed, let us consider the Fourier
transform of (\ref{eq:1.1}). For a Fourier mode $\hat{a}_j(k,t)$ we
obtain
$$ \frac{d}{dt} (\phi\, R_j \, \hat{a}_j)
= - \phi\, (\lambda_j \, R_j + k^2 \, D_c - i k \, \U) \, \hat{a}_j
+ y_i \, R_{j-1} \, \phi \,  \hat{a}_{j-1} = - \lambda'_j \, R_j
 \, \phi\, \hat{a}_j + y_i \, R_{j-1} \phi\,  \hat{a}_{j-1}
$$
with $\lambda'_j = \lambda_j + (k^2 \, D_c -ik  \, \U )/R_j$ for
$j=1,..,N_r-1$. A  transform analogous to (\ref{eq:1.9}) uncouples
the problem. By taking the partial differential equation
counterpart, we obtain an equation similar to  (\ref{eq:1.10}).

\end{rem}

\section{Discrete tools}
\label{sec:dtools}

\newcounter{disc}
\renewcommand{\theequation}{\thesection.\thedisc}

\subsection{Mesh and discrete spaces}
\label{subsec:maillage}

Let ${\cal{T}}_h$ be a triangular mesh of $\Omega$. The
circumscribed circle of a triangle $K \in {\mathcal{T}}_h$ is
centered at  $\x_K$ and has the diameter $h_K$. We set $h=\max_{K
\in {\mathcal{T}}_h} h_K$. We assume that all the interior angles of
the triangles of the mesh are less than $\frac{\pi}{2}$, so that
$\x_K \in K$. The set of the edges of the triangle $K \in
{\mathcal{T}}_h$ is ${\cal{E}}_K$. The symbol $\N_{K,\sigma}$
denotes the unit normal vector  to an edge $\sigma \in {\cal{E}}_K$
and pointing outward $K$. We denote by ${\cal{E}}_h$ the set of the
edges of the mesh.
 We distinguish the subset ${\mathcal{E}}^{int}_h \subset {\mathcal{E}}_h$ (resp. ${\mathcal{E}}^{ext}_h$)
 of the edges located inside  $\Omega$ (resp. on $\partial \Omega$).
The middle of an  edge $\sigma \in {\cal{E}}_h$ is $\x_\sigma$ and
its length is $|\sigma|$. For each edge $\sigma \in
{\cal{E}}^{int}_h$
 let $K_\sigma$ and $L_\sigma$ be the two  triangles having $\sigma$ in
 common; we set
 $d_\sigma=d(\x_{K_\sigma},\x_{L_\sigma})$. For all $\sigma \in {\mathcal{E}}^{ext}_h$ only the triangle $K_\sigma$
  located inside $\Omega$ is defined and we set $d_\sigma=d(\x_{K_\sigma},\x_{\sigma})$. Then for all $\sigma \in {\mathcal{E}}_h$ we set
   $\tau_\sigma=\frac{|\sigma|}{d_\sigma}$.
We assume the following  on the mesh (see \cite{eym}): there exists
$C>0$ such that
$$
   \forall \, \sigma \in {\cal{E}}_h \,,  \hspace{1cm}  d_\sigma \ge C \, |\sigma|
  \quad \mathrm{and} \quad
 |\sigma| \ge C \, h.
$$
It implies that there exists $C>0$ such that \stepcounter{disc}
 \begin{equation}
 \label{eq:mintaus}
  \forall \, \sigma \in {\cal{E}}^{int}_h \, , \hspace{1cm} \tau_\sigma=|\sigma|/d_\sigma \ge C.
 \end{equation}
We define on the mesh the following spaces. The usual space for
finite volume schemes is
 \begin{equation*}
  P_0 = \{ q \in L^2 \; ; \; \forall \, K \in {\cal{T}}_h, \; \; q|_K \hbox{ is a  constant} \}.
\end{equation*}
For any function $q_h \in P_0$ and any  $K \in {\cal{T}}_h$ we set $
q_K=q_h|_K$. We also consider
\begin{eqnarray*}
  P^d_1 &=& \{ q \in L^2 \; ; \; \forall \, K \in {\cal{T}}_h, \; \; q|_K \hbox{ is affine}
  \} \, , \\
  P^c_1 &=& \{ q_h \in P^d_1 \; ; \; q_h \hbox{ is continuous over }\Omega   \}\,   ,\\
  P^{nc}_1 &=& \{ q_h \in P^d_1 \; ; \; \forall \, \sigma \in {\cal{E}}^{int}_h, \, q_h \hbox{ is continuous at the middle of  } \sigma \}.
 \end{eqnarray*}
 We have $P_1^c \subset H^1$. On the other hand $P^{nc}_1 \not \subset H^1$, but $P^{nc}_1 \subset H_d^1$ with
\begin{equation*}
  H^1_d= \{ q \in L^2 \; ; \; \forall \, K \in {\cal{T}}_h, \; \; q|_K \in H^1(K) \}.
\end{equation*}
Thus we define
  $\widetilde \nabla_h:   H^1_d \to \L^2$ by setting
\stepcounter{disc}
\begin{equation}
\label{eq:defgradht}
  \forall \, q_h \in  H_d^1, \quad  \forall \, K \in {\cal{T}}_h, \quad \widetilde \nabla_h   q_h|_K = \nabla q_h|_K
\end{equation}
and the associated  norm $\|.\|_{1,h}$  is given by
$$
\forall \, q_h \in H_d^1, \qquad \|q_h\|_{1,h}=( |q_h|^2 +
|\widetilde \nabla_h q_h|^2 )^{1/2}.
$$
We then have the  following Poincar\'e-like inequality for the space
$P^{nc}_1 \cap L^2_0$ (see \cite{ach}).
\begin{prop}
\label{proppoinp1nc} There exists  $C>0$ such that $|q_h| \le C \,
|\widetilde \nabla_h q_h|$ for all $q_h \in P^{nc}_1 \cap L^2_0$.
\end{prop}
\noindent We also define discrete analogues of the norms $H^1$ and
$H^{-1}$ for the space $P_0$ by setting
$$
  \|q_h\|_h=\Bigl( \sum_{\sigma \in {\cal{E}}^{int}_h} \tau_{\sigma} \, (q_{L_\sigma} - q_{K_\sigma})^2 \Bigr)^{1/2}
\quad \mathrm{and} \quad
  \|q_h\|_{-1,h}=\sup_{\psi_h \in P_0} \frac{(q_h,\psi_h)}{\|\psi_h\|_h}
$$
for any function $q_h \in P_0$. Note that for any $p_h \in P_0$ and
$q_h \in P_0$, $  (p_h,q_h)  \le  \|q_h\|_{-1,h} \, \|q_h\|_h$. The
following Poincar\'e-like inequality holds for the space $P_0 \cap
L^2_0$ (see \cite{eym}).
\begin{prop}
\label{proppoinp0} There exists $C>0$ such that $|q_h| \le C \,
\|q_h\|_h$ for all $q_h \in P_0 \cap L^2_0$.
\end{prop}
\noindent Finally we set $\P_0=(P_0)^2$, $\P_1^d=(P_1^d)^2$ and use
the Raviart-Thomas spaces
 \cite{brezzi}
\begin{align*}
  \mathbf{RT^d_0}&= \{ \V_h \in \P^d_1 \; ; \quad \forall \, K \in {\mathcal{T}}_h \, , \quad \forall \, \sigma \in {\cal{E}}_K, \quad \V_h|_K \cdot \N_{K,\sigma}
   \hbox{ is a  constant}  \} \, , \\
  \mathbf{RT_0}&= \{ \V_h \in \mathbf{RT^d_0} \; ; \quad  \forall \, \sigma
   \in {\cal{E}}^{int}_h,
 \quad   \V_h|_{K_\sigma} \cdot \N_{K_\sigma,\sigma} = \V_h|_{L_\sigma} \cdot \N_{K_\sigma,\sigma} \quad \hbox{ and } \quad
  \V_h \cdot \N|_{\partial \Omega}=0\}.
\end{align*}
For all $\V_h \in \mathbf{RT_0}$,  $K \in {\cal{T}}_h$ and  $\sigma
\in {\cal{E}}_K$, we set $  (\V_h \cdot
\N_{K,\sigma})_\sigma=\V_h|_K \cdot \N_{K,\sigma}$.

\subsection{Projection operators}
\label{subsec:opproj}

We associate with the spaces of section \ref{subsec:maillage} some
projection operators. First, we define  $\Pi_{P_1^c}: H^1_d \to
P_1^c$
by setting
 \stepcounter{disc}
\begin{equation}
\label{eq:defprojp1c} \forall \, q \in H^1_d \, , \hspace{1cm}
\forall \, \phi_h \in P^c_1 \, , \hspace{1cm} \big(\nabla
(\Pi_{P^c_1} q),\nabla \phi_h\big)=(\nabla q, \nabla \phi_h).
\end{equation}
Next, we consider the space $P_0$. Let ${\mathcal{C}}_d=\{q_h \in
L^2 \, ; \hspace{.3cm} q_h \hbox{ is equal a.e. to a continuous
function}\}$. We define   $\Pi_{P_0}:  L^2 \to P_0$ and $\widetilde
\Pi_{P_0}:{\mathcal{C}}_d \to P_0$ by setting \stepcounter{disc}
\begin{equation}
\label{eq:defprojp0}
 (\Pi_{P_0} p)_K = \frac{1}{|K|} \int_K p(\x) \, d\x \, ,
\hspace{1cm}  (\widetilde \Pi_{P_0} q)_K = q(\x_K) \, ,
\end{equation}
for all $p \in L^2$, $q \in {\mathcal{C}}_d$ and $K \in
{\mathcal{T}}_h$. We also set $\Pi_{\P_0}=(\Pi_{P_0})^2$.
 For the space $P^{nc}_1$, we define $\widetilde \Pi_{P^{nc}_1}:L^2 \to P^{nc}_1$ and $\Pi_{P^{nc}_1}:H^1 \to P^{nc}_1$. For all
 $p \in L^2$ and $q \in H^1$,  $\widetilde \Pi_{P^{nc}_1} p$ and $\Pi_{P^{nc}_1} q$ satisfy
\stepcounter{disc}
\begin{equation}
\label{eq:defprojp1nc} \forall \, \psi_h \in P^{nc}_1 \, , \quad
(\widetilde \Pi_{P^{nc}_1} p,\psi_h)=(p,\psi_h) \, ; \hspace{1cm}
\forall \, \sigma \in {\mathcal{E}}_h \, , \quad \int_\sigma
(\Pi_{P^{nc}_1} q) \, d\sigma   =\int_\sigma q \, d\sigma.
\end{equation}
For the space $\mathbf{RT_0}$, we define $\widetilde
\Pi_{\mathbf{RT_0}}:\L^2 \to \mathbf{RT_0}$ and
$\Pi_{\mathbf{RT_0}}:\H^1 \to \mathbf{RT_0}$. For all
 $\V \in \L^2$ and  $\W \in \H^1$, $\widetilde \Pi_{\mathbf{RT_0}} \V$ and $\Pi_{\mathbf{RT_0}} \W$ satisfy
\begin{equation}
\stepcounter{disc} \label{eq:defprojrt0} \forall \, \W_h \in
\mathbf{RT_0} \, , \quad (\widetilde \Pi_{\mathbf{RT_0}}
\V,\W_h)=(\V,\W_h) \, ; \hspace{.2cm}   \forall \, \sigma \in
{\mathcal{E}}^{int}_h \, , \quad  \int_\sigma
(\W-\Pi_{\mathbf{RT_0}} \W) \cdot \N_{K_\sigma,\sigma} \, d\sigma
   =0.
\end{equation}
The operators $\Pi_{P_0}$, $\widetilde \Pi_{P^{nc}_1}$ (resp.
$\Pi_{\P_0}$, $\widetilde \Pi_{\mathbf{RT_0}}$) are $L^2$ (resp.
$\L^2$) projection operators. They are stable  for the $L^2$ (resp.
$\L^2$) norms. The operators $\widetilde \Pi_{P_0}$,
$\Pi_{P^{nc}_1}$ and $\Pi_{\mathbf{RT_0}}$ are interpolation
operators.
The following estimates are classical (\cite{brenn} p.109 and
\cite{brezzi}).
\begin{prop}
\label{prop:errint} 
There exists $C>0$ such that for all  $q \in H^1$ and $\V \in \H^1$
\begin{equation*}
 |q-\Pi_{P_0} q|\le C \, h \, \|q\|_1  , \hspace{1cm}
   |\V-\Pi_{\mathbf{RT_0}} \V| \le C \, h \, \|\V\|_1.
\end{equation*}
For all $p \in H^1$ and $q \in H^2$ we have
\begin{equation*}
 |p-\Pi_{P^{nc}_1} p|\le C \, h \, \|p\|_1 \, , \hspace{1cm}
  |\widetilde \nabla_h(q-\Pi_{P^{nc}_1} q)| \le C \, h \, \|q\|_2.
\end{equation*}
For all $q \in H^d_1$ we have
\begin{equation*}
  |q-\Pi_{P^{c}_1} q| \le C \, h \, \|q\|_{1,h}.
\end{equation*}
%
\end{prop}
\noindent Finally, using the Sobolev embedding theorem, one  checks
that \stepcounter{disc}
\begin{equation}
\label{eq:esti}
  |\Pi_{P_0}  q - \widetilde \Pi_{P_0} q | \le C \,  h \, \|q\|_{W^{1,r}}
\end{equation}
for all $q \in W^{1,r}$ with $r>1$ (see \cite{zimconv}).

\subsection{Discrete operators}
\label{subsec:defop}

Equations (\ref{eq:1.10})--(\ref{eq:1.12}) use the differential
operators gradient, divergence and laplacian. We have to define
analogous operators in the discrete  setting.
 The discrete gradient operator $\nabla_h:P^{nc}_1 \to P_0$ is the restriction to $P^{nc}_1$ of the operator $\widetilde \nabla_h$ given by (\ref{eq:defgradht}).
The discrete divergence operator $ \hbox{div}_h: \P_0 \to P^{nc}_1$
is defined by
\begin{equation*}
\begin{array}{lcl}
&& {\displaystyle
 \forall  \sigma \in {\cal{E}}^{int}_h, \quad
(\hbox{div}_h  \, \V_h) (\x_\sigma)= \frac{3 \,
|\sigma|}{|K_\sigma|+|L_\sigma|} \, (\V_{L_\sigma} - \V_{K_\sigma})
\cdot \N_{K,\sigma} , }
 \\
&&   {\displaystyle \forall  \sigma \in {\cal{E}}^{ext}_h, \quad
(\hbox{div}_h  \, \V_h) (\x_\sigma)= -\frac{3 \,
|\sigma|}{|K_\sigma|} \, \V_{K_\sigma} \cdot \N_{K,\sigma}, }
\end{array}
\end{equation*}
for all $\V_h \in \P_0$. It is adjoint to $\nabla_h$ (proposition
\ref{prop:propadjh} below).
 The discrete laplacian operator $ \Delta_{h}: P_0 \to P_0$ is the usual one for finite volume schemes (see \cite{eym}). For all
 $q_h \in P_0$ and $K \in {\mathcal{T}}_h$ we have
\stepcounter{disc}
\begin{equation}
\label{eq:deflap}
  \Delta_{h}  q_h  \big|_K=
   \frac{1}{|K|}\sum_{\sigma \in {\cal{E}}_K \cap {\cal{E}}^{int}_h} \tau_\sigma  \, \big( q_{L_\sigma}-q_{K_\sigma}
  \big).
\end{equation}
Let us now consider the convection terms in (\ref{eq:1.10}) and
(\ref{eq:1.11}). We define  $\bt: \H^1 \times H^1 \to L^2$ by
 \stepcounter{disc}
\begin{equation}
\label{eq:defbt} \bt(\V,q)= \hbox{div}(q  \, \V)
\end{equation}
for all $q \in H^1$ and $\V \in \H^1$. In order to define a discrete
counterpart to $\bt$ we use the classical upwind scheme (see
\cite{eym}). The discrete operator $\bt_h: \mathbf{RT_0} \times P_0
\to P_0$ is such that
 \stepcounter{disc}
\begin{equation}
\label{eq:defbth}
  \bt_h(\V_h,q_h)\big|_K=
\frac{1}{|K|} \sum_{\sigma \in {\cal{E}}_K \cap {\cal{E}}^{int}_h}
\,
   |\sigma| \, \bigl( (\V_h \cdot \N_{K,\sigma})^+_\sigma \, q_K +
(\V_h \cdot \N_{K,\sigma})^-_\sigma \, q_{L_\sigma}  \bigr)
\end{equation}
for all $\V_h \in \mathbf{RT_0}$, $q_h \in P_0$ and $K \in
{\mathcal{T}}_h$.
We have set $a^+=\max(a,0)$ and $a^-=\min(a,0)$ for all $a \in
\mathbb{R}$. Integrating by parts the convection terms also leads to
consider  $\b: \L^2 \times L^2 \times L^\infty \to \mathbb{R}$
defined by
 \stepcounter{disc}
\begin{equation}
\label{eq:defb} \b(\V,p,q) = - \int_\Omega p \, \V \cdot \nabla q \,
d\x
\end{equation}
for all $\V \in L^2$, $p \in L^2$ and $q \in L^\infty$. The discrete
counterpart is $ \b_h: \mathbf{RT_0} \times P_0 \times P_0  \to
\mathbb{R}$ with
 \stepcounter{disc}
\begin{equation}
\label{eq:defbh}
  \b_h(\V_h,p_h,q_h)=\sum_{K \in {\cal{T}}_h} q_K
  \sum_{\sigma \in {\cal{E}}_K \cap {\cal{E}}^{int}_h} \,
   |\sigma| \, \bigl( (\V_h \cdot \N_{K,\sigma})^+_\sigma \, p_K +
(\V_h \cdot \N_{K,\sigma})^-_\sigma \, p_{L_\sigma}  \bigr)
\end{equation}
for all $\V_h \in \mathbf{RT_0}$, $p_h \in P_0$ and $q_h \in P_0$.


\section{Properties of the discrete operators}
\label{sec:propop}
\newcounter{propop}
\renewcommand{\theequation}{\thesection.\thepropop}

 The properties of the discrete operators are analogous to the ones satisfied by their continuous counterpart.
The gradient and divergence operators are adjoint. For the operators
$\nabla_h$ and $\hbox{div}_h$ we
  state in \cite{zimstab} the following.
\begin{prop}
\label{prop:propadjh} For all $\V_h \in \P_0$ and $q_h \in P^{nc}_1$
we have
$  (\V_h, \nabla_h q_h)=-(q_h,\hbox{\rm div}_h  \V_h)$.
\end{prop}

\noindent Let us now consider the convection terms. Let  $q \in
L^\infty \cap H^1$, $\V \in \L^2$ with $\hbox{div} \, \V \in L^2$
and  $\hbox{div} \, \V(\x) \ge 0$ a.e. in $\Omega$. We obtain $
\b(\V,q,q)=\int_\Omega  (q^2/2)  \, \hbox{div}  \,  \V  \, d\x \ge
0$ by integration by parts. For $\b_h$ we state in \cite{zimstab} a
similar result.

\begin{prop}
\label{prop:posbh} Let $\V_h \in \mathbf{RT_0}$ with $\hbox{\rm div}
\, \V_h \ge 0$. We have $\b_h(\V_h,q_h,q_h) \ge 0$ for all $q_h \in
P_0$
.
\end{prop}

\noindent The following stability properties are used to prove the
error estimates in section \ref{sec:esterr}.
\begin{prop}
\label{prop:stabbth} There exists $C>0$ such that for all  $p_h \in
P_0$, $q_h \in P_0$ and $\V_h \in \mathbf{RT_0}$ with  $\hbox{ \rm
div} \,   \V_h=0$ \stepcounter{propop}
\begin{equation}
\label{eq:majbh} |\b_h(\V_h,p_h,q_h)| \le C \, |\V_h| \, \|p_h\|_{h}
\, \|q_h\|_h.
\end{equation}
There exists $C>0$ such that for all  $p_h \in P_0$, $q_h \in P_0
\cap L^2_0$, $\V_h \in \mathbf{RT_0}$ \stepcounter{propop}
\begin{equation}
\label{eq:majbhb} |\b_h(\V_h,p_h,q_h)| \le C \, (|\V_h| \,
\|p_h\|_h+|\hbox{\rm div}\,   \V_h| \, \|p_h\|_{\infty}\, ) \,
\|q_h\|_{h}.
\end{equation}
\end{prop}
\noindent {\sc Proof.} For all $K \in {\cal{T}}_h$ and $\sigma \in
{\cal{E}}_K \cap {\cal{E}}^{int}_h$ we write
\begin{equation*}
  (\V_h \cdot \N_{K,\sigma})^+_\sigma \, p_K +  (\V_h \cdot \N_{K,\sigma})^-_\sigma
  \, p_{L_\sigma}
  =(\V_h \cdot \N_{K,\sigma})_\sigma \, p_K
 - |(\V_h \cdot \N_{K,\sigma})_\sigma| \, (p_{L_\sigma}-p_K).
\end{equation*}
Thus  (\ref{eq:defbh}) reads $\b_h(\V_h,p_h,q_h)=S_1+S_2$ with
\begin{equation*}
  S_1 = -\sum_{K \in {\cal{T}}_h} q_K \sum_{\sigma \in {\cal{E}}_K
\cap {\cal{E}}^{int}_h} |\sigma|\, |(\V_h \cdot
\N_{K,\sigma})_\sigma| \, (p_{L_\sigma}-p_K) \, ,  \hspace{.7cm}
S_2= \sum_{K \in {\cal{T}}_h} p_K \, q_K \sum_{\sigma \in
{\cal{E}}_K \cap {\cal{E}}^{int}_h} |\sigma| \, (\V_h \cdot
\N_{K,\sigma})_\sigma. \label{eq:decompbh}
\end{equation*}
Using  the Cauchy-Schwarz inequality we write
\begin{eqnarray*}
&& \vert S_1 \vert= \Bigl\vert \sum_{\sigma \in  {\cal{E}}^{int}_h}
  |\sigma|\, |\V_h(\x_\sigma) \cdot
\N_{K,\sigma}| \, (p_{L_\sigma}-p_K) \, (q_{L_\sigma}-q_K)
\Bigr\vert
\\
&& \qquad
 \le  h \, \|\V_h\|_{\infty} \, \bigl( \sum_{\sigma \in {\cal{E}}^{int}_h} (p_{L_\sigma} - p_{K_\sigma})^2 \bigr)^{1/2}
  \bigl( \sum_{\sigma \in {\cal{E}}^{int}_h} (q_{L_\sigma} - q_{K_\sigma})^2 \bigr)^{1/2}.
\end{eqnarray*}
Since $\V_h \in \mathbf{RT_0} \subset (P^d_1)^2$  we have   $h \,
\|\V_h\|_{\infty} \le C \, |\V_h|$ (\cite{brenn} p. 112). Moreover
(\ref{eq:mintaus}) implies $
 \sum_{\sigma \in {\cal{E}}^{int}_h} (p_{L_\sigma} - p_{K_\sigma})^2
 \le  C  \, \sum_{\sigma \in {\cal{E}}^{int}_h} \tau_\sigma \, (p_{L_\sigma} - p_{K_\sigma})^2
 = C  \,  \|p_h\|^2_h
$ and $
 \sum_{\sigma \in {\cal{E}}^{int}_h} (q_{L_\sigma} - q_{K_\sigma})^2
 \le C \, \|q_h\|^2_h
$. Thus \stepcounter{propop}
\begin{equation}
\label{eq:estbh6}
  |S_1| \le C \, |\V_h| \, \|p_h\|_h \, \|q_h\|_h.
\end{equation}
We now consider $S_2$. We have $(\V_h \cdot \N_{K,\sigma})_\sigma=0$
for all $K \in {\cal{T}}_h$ and $\sigma \in {\cal{E}}_K \cap
{\cal{E}}^{ext}_h$. Thus we write
\begin{equation*}
\sum_{\sigma \in {\cal{E}}_K \cap {\cal{E}}^{int}_h} |\sigma| \,
(\V_h \cdot \N_{K,\sigma})_\sigma =  \sum_{\sigma \in {\cal{E}}_K}
|\sigma| \, (\V_h \cdot \N_{K,\sigma})_\sigma =\int_K \hbox{div}  \,
\V_h \, d\x.
\end{equation*}
It gives the following relation.
\begin{equation*}
\label{eq:estbh1} S_2 = \sum_{K \in {\cal{T}}_h} p_K \, q_K \,
\int_K \hbox{div}  \V_h \, d\x=\int_\Omega p_h \, q_h \, \hbox{div}
\V_h  \, d\x.
\end{equation*}
Thus if $\hbox{div} \,  \V_h=0$ then  $S_2=0$ and estimate
(\ref{eq:estbh6}) gives (\ref{eq:majbh}). Let us prove
(\ref{eq:majbhb}).
Since $q_h \in P_0 \cap L^2_0$ we can apply proposition \ref{proppoinp0}. Using the Cauchy-Schwarz inequality we get
\stepcounter{disc}
\begin{equation*}
  |S_2| \le \|p_h\|_\infty \, |q_h| \, |\hbox{div}  \, \V_h|
  \le C \, \|p_h\|_\infty \, \|q_h\|_h \, |\hbox{div} \,  \V_h|.
  \label{eq:estbh66}
\end{equation*}
This latter estimate together with (\ref{eq:estbh6}) gives
(\ref{eq:majbhb}).
 \qed

\noindent Lastly, we claim that $\bt_h$ is a consistent
approximation of $\bt$ \cite{zimstab}.
\begin{prop}
\label{propconsbh} Let $r>0$. There exists $C>0$ such that for all
functions $q \in H^2$ and  $\V \in \H^{1+r}$ with $\V \cdot
\N|_{\partial \Omega}=0$
\begin{equation*}
\label{eq:consbh}
  \|\Pi_{P_0}\bt(\V,q) - \bt_h(\Pi_{\mathbf{RT_0}} \V, \widetilde \Pi_{P_0} q)\|_{-1,h}
\le C \, h  \,  \|q\|_1 \,  \|\V\|_{1+r}.
\end{equation*}
\end{prop}

\noindent Let us now consider the discrete laplacian. We have a
coercivity  and stability result.
\begin{prop}
\label{prop:coerlap} For all $p_h \in P_0$ and $q_h \in P_0$, we
have
\begin{equation*}
 -( \Delta_{h} p_h,p_h)=  \|p_h\|^2_h \, , \hspace{1cm}
\vert (\Delta_{h} p_h,q_h) \vert \le \|p_h\|_h \, \|q_h\|_h .
\end{equation*}
\end{prop}
\noindent {\sc Proof.} Definition (\ref{eq:deflap}) implies
\stepcounter{propop}
\begin{equation}
  (\Delta_{h} p_h,q_h) =
  \sum_{K \in {\cal{T}}_h} \! q_K  \!\! \sum_{\sigma \in {\cal{E}}_K \cap {\cal{E}}^{int}_h} \!\!\! \tau_\sigma  \,  (p_{L_\sigma}-p_{K_\sigma})
= - \sum_{\sigma \in {\cal{E}}^{int}_h} \tau_\sigma \,
(p_{L_\sigma}-p_{K_\sigma})  \, (q_{L_\sigma}-q_{K_\sigma}).
\label{eq:propdh3}
\end{equation}
Setting $q_h=p_h$ gives the first part of the result. Using the
Cauchy-Schwarz inequality, we get the second one.  \qed

\noindent We also deduce from (\ref{eq:propdh3}) the following
property.
\begin{prop}
\label{prop:adjlap} For all $p_h \in P_0$ and $q_h \in P_0$ we have
  $(\Delta_{h} p_h, q_h)=(p_h,\Delta_{h} q_h)$.
\end{prop}

\noindent Lastly, we state that $\Delta_{h}$ is a consistent
approximation of the laplacian. The proof follows the lines of the
one of proposition 1.14 in \cite{zimconv}.
\begin{prop}
\label{propconslap} There exists $C>0$ such that for all $q \in H^2$
with  $\nabla q \cdot \N|_{\partial \Omega}=0$ we have
\stepcounter{disc}
\begin{equation*}
\label{eq:conslap}
 \|\Pi_{P_0} (\Delta q) - \Delta_{h}(\widetilde \Pi_{P_0} q)\|_{-1,h} \le C \,  h \,
\|q\|_2.
\end{equation*}
\end{prop}

\section{The finite element-finite volume scheme}
\label{sec:schema}

\newcounter{sch}
\renewcommand{\theequation}{\thesection.\thesch}

We now introduce the  scheme for (\ref{eq:1.10})-(\ref{eq:1.14}).
The interval $[0,T]$ is split with a constant time step  $k=T/N$. We
set $ [0,T]=\bigcup^{N-1}_{m=0} [t_m,t_{m+1}]$ with $t_m =m \, k$.
The time derivatives are approximated using a first order Euler
scheme. The convection terms are discretized semi-implicitly in time
and the other ones in an implicit way.
We set $s_h = \Pi_{P_0} s $, $ s_h^{c} = \Pi_{P_0}
s_{c}$, $s_h^{\theta} = \Pi_{P_0} s_\theta$ and  $\F^m_h=\Pi_{\P_0}
\F(t_m)$ for all $m \in \{0, \dots, N \}$. Since $\Pi_{P_0}$ (resp.
$\Pi_{\P_0}$) is stable for the $L^2$ (resp. $\L^2$) norm we have
\stepcounter{sch}
\begin{equation}
\label{eq:stabfh}
 |s_h| \le |s| \, , \quad |s^c_h| \le |s_c| \, , \quad |s^\theta_h| \le |s_\theta| \, , \quad|\F^m_h| \le |\F(t_m)| \le \|\F\|_{L^\infty(0,T;\L^2)}.
\end{equation}
The initial values are  $c^0_h=\Pi_{P_0} c_{0}$ and
$\theta^0_h=\Pi_{P_0} \theta_0$.
Then for all $n \in \{ 0, \dots, N-1 \}$, the quantities
$c^{n+1}_{h}\in P_0 $, $\theta^{n+1}_h \in P_0$, $p^{n+1}_h \in
P^{nc}_1 \cap L^2_0$, $\U^{n+1}_h \in \mathbf{RT_0} $ are the
solutions of the following problem.\stepcounter{sch}
\begin{eqnarray}
&&  \frac{c^{n+1}_{h}-c^n_h}{k} -D_c  \, \Delta_{h} c^{n+1}_h
=s_h^{c} -(s_h+\lambda) \, c^{n+1}_h -
  \bt_h(\U^n_h,c^{n+1}_h)  ,
\label{eq:eqcd} \\
\stepcounter{sch} && \frac{\theta^{n+1}_{h}-\theta^n_h}{k} -D_\theta
\, \Delta_{h} \theta^{n+1}_h =- s_h^{\theta} - s_h \,  (
\theta^{n+1}_h-\theta_*) -
  \bt_h(\U^n_h,\theta^{n+1}_h) ,
\label{eq:eqtd} \\
\stepcounter{sch} &&  \hbox{div}_h ( \kappa^{n+1}_h  \, \nabla_h
p^{n+1}_h )=
  \hbox{div}_h \, \F^{n+1}_h - \widetilde \Pi_{P^{nc}_1} s_h \,  ,
\label{eq:eqpd} \\
\stepcounter{sch} &&  \U^{n+1}_h = \widetilde \Pi_{\mathbf{RT_0}} (
\F^{n+1}_h - \kappa^{n+1}_h  \, \nabla_h p^{n+1}_h ),
\label{eq:equd}
\end{eqnarray}
with $\kappa^{n+1}_h=\kappa(c^{n+1}_h, \theta^{n+1}_h) \in P_0$.
This term is defined thanks to proposition \ref{propmaxd} below.
Note also that the boundary conditions are implicitly included in
the definition of the discrete operators (section
\ref{subsec:defop}). The existence of a unique solution to
(\ref{eq:eqcd}) and (\ref{eq:eqtd}) is classical (see \cite{eym}).
Since $\kappa^{n+1}_h \ge \kappa_{min}>0$ and $p^{n+1}_h \in L^2_0$
equation (\ref{eq:eqpd}) also has a unique solution (see
\cite{brenn}). We have a discrete equivalent for the divergence
condition (\ref{eq:1.12}).
\begin{prop}
\label{prop:propdivsh}
 For all $m \in \{ 1, \dots, N \}$ we have 
$ \hbox{\rm div} \,  \U^m_h = s_h$ .
\end{prop}
%
\noindent {\sc Proof.}  Let $m\in\{1,\dots,N\}$ and $n=m-1$. We
compare the solution of (\ref{eq:eqpd})--(\ref{eq:equd}) with the
solution of the following mixed hybrid problem. Let
 ${\cal{E}}_0=\left\{\mu_h: \cup_{\sigma \in {\cal{E}}_h} \to \mathbb{R} \, ; \; \forall \, \sigma \in {\cal{E}}_h \, ,
 \; \mu_h|_\sigma \hbox{ is constant} \right\}$.
Then $\Ut^{m}_h \in \mathbf{RT^d_0}$, $\overline{p}^m_h \in P_0$ and
$\lambda^m_h \in {\cal{E}}_0$ are the solution of (see
\cite{brezzi}) \stepcounter{sch}
\begin{eqnarray}
\hspace{-1.2cm} &&  \forall  \, \V_h \in \mathbf{RT^d_0} \,,
\hspace{.2cm} (\Ut^m_h,\V_h) + \sum_{K \in {\cal{T}}_h} \kappa^m_K
\sum_{\sigma' \in {\cal{E}}_K} |\sigma'| \, \lambda^m_{\sigma'} \,
(\V_h|_K \cdot \N_{K,\sigma'}) - \sum_{K \in {\cal{T}}_h} |K|
\kappa^m_K \, \overline{p}^m_K \, \hbox{div}  \, \V_h|_K
 = (\F^m_h,\V_h) \, ,
\label{eq:mha} \\
\stepcounter{sch} \hspace{-1.2cm}&&    \forall  \, \mu_h \in
{\cal{E}}_0, \qquad \sum_{K \in {\cal{T}}_h} \int_{\partial K} \mu_h
\,  (\Ut^m_h \cdot \N) \, d\sigma =0  \, , \hspace{1cm}  \forall  K
\in {\cal{T}}_h \, , \quad
  \int_K \hbox{div}  \, \Ut^m_h \, d\x=
\int_K s \, d\x \, ,
 \label{eq:mhc}
\end{eqnarray}
and $\widetilde p^m_h \in P^{nc}_1$ is defined by
 $\int_\sigma \widetilde p^m_h \, d\sigma = \lambda^m_\sigma$
for all  $\sigma \in {\cal{E}}_h$.
 Let  $\sigma \in {\cal{E}}_h$. We define    $\phi_\sigma \in
P^{nc}_1$ by setting $\phi_\sigma(\x_\sigma)=1$ and
$\phi_\sigma(\x_{\sigma'})=0$ for all $\sigma' \in
{\cal{E}}_h\backslash\{\sigma\}$. We set   $\V_h=\nabla_h
\phi_\sigma \in \P_0 \subset \mathbf{RT^d_0}$ in (\ref{eq:mha}). We
have
\begin{equation*}
 \sum_{K \in {\cal{T}}_h} \kappa^m_K \sum_{\sigma' \in {\cal{E}}_K}
|\sigma'| \, \lambda^m_{\sigma'} \, \nabla_h \phi_\sigma|_K \cdot
\N_{K,\sigma'} = \sum_{K \in {\cal{T}}_h} \kappa^m_K \, \nabla_h
\phi_\sigma|_K \cdot \sum_{\sigma' \in {\cal{E}}_K} |\sigma'| \,
\lambda^m_{\sigma'} \, \N_{K,\sigma'}
\end{equation*}
and according to the gradient formula
\begin{equation*}
  \sum_{\sigma' \in {\cal{E}}_K} |\sigma'| \,
\lambda^m_{\sigma'} \, \N_{K,\sigma'}
  =\sum_{\sigma' \in {\cal{E}}_K} \int_{\sigma'}
 \widetilde p^m_h \, \N_{K,\sigma'} \, d{\sigma'}
 = \int_K \nabla_h  \widetilde
 p^m_h\, d\x.
\end{equation*}
Thus we get from (\ref{eq:mha})
  \stepcounter{sch}
 \begin{equation}
 \label{eq:lienmh1}
(\Ut^m_h,\nabla_h \phi_\sigma) +(\kappa^m_h \, \nabla_h \widetilde
p^m_h , \nabla_h \phi_\sigma) =(\F^{m}_h,\nabla_h \phi_\sigma).
 \end{equation}
The first term  in (\ref{eq:lienmh1}) is treated as follows.
Integrating by parts we get
\begin{equation*}
(\Ut^m_h,\nabla_h \phi_\sigma) =-(\phi_\sigma, \hbox{div} \,
\Ut^m_h) + \sum_{K \in {\cal{T}}_h} \sum_{\sigma' \in {\cal{E}}_K}
\int_{\sigma'} \phi_\sigma \, (\Ut^m_h|_K \cdot \N_{K,\sigma'}) \,
d\sigma'.
\end{equation*}
Since (\ref{eq:mhc}) implies that $\Ut^m_h \in \mathbf{RT_0}$, we
have
\begin{equation*}
\sum_{K \in {\cal{T}}_h} \sum_{\sigma' \in {\cal{E}}_K}
\int_{\sigma'} \phi_\sigma \, (\Ut^m_h \cdot \N_{K,\sigma'}) \,
d\sigma' =\sum_{\sigma \in {\cal{E}}^{int}_h} |\sigma| \,
\phi_\sigma(\x_\sigma) \, (\Ut^m_h|_{L_\sigma}\cdot
\N_{K_\sigma,\sigma} - \Ut^m_h |_{K_\sigma} \cdot
\N_{K_\sigma,\sigma})=0.
\end{equation*}
Thus $(\Ut^m_h,\nabla_h \phi_\sigma)=-(\phi_\sigma,\hbox{div} \,
\Ut^m_h)$. Then, using (\ref{eq:mhc}), we get \stepcounter{sch}
\begin{equation*}
\label{eq:lienmh4} (\Ut^m_h,\nabla_h \phi_\sigma)
=-(\phi_\sigma,s_h) =-(\phi_\sigma,\widetilde \Pi_{P^{nc}_1} s_h).
\end{equation*}
Furthermore, according to proposition \ref{prop:propadjh}, we have
$(\kappa^m_h \nabla_h \widetilde p^m_h , \nabla_h \phi_\sigma)
=-\big(\phi_\sigma, \hbox{div}_h  (\kappa^m_h \, \nabla_h \widetilde
p^m_h)\big)$ and $(\F^{m}_h,\nabla_h
\phi_\sigma)=-(\phi_\sigma,\hbox{div}_h  \F^m_h)$. Hence we deduce
from (\ref{eq:lienmh1}) that
\begin{equation*}
\forall \phi_\sigma \in P^{nc}_1 \, , \hspace{1cm} \big(
\phi_\sigma,\hbox{div}_h (\kappa^m_h \nabla_h \widetilde
p^m_h)-\hbox{div}_h  \F^m_h + \widetilde \Pi_{P^{nc}_1}
s_h^m\big)=0.
\end{equation*}
Since $(\phi_\sigma)_{\sigma \in {\cal{E}}_h}$ is a basis of
$P^{nc}_1$, we get $\hbox{div}_h (\kappa^m_h \nabla_h \widetilde
p^m_h)=\hbox{div}_h \F^m_h -\widetilde\Pi_{P^{nc}_1} s_h^m$. Thus,
by (\ref{eq:eqpd}), there exists a real $C$ such that
 $\widetilde p^m_h = p^m_h +C$.
We now compare $\Ut^m_h$ with $\U^m_h$. Since for all  $\V_h \in
\mathbf{RT_0}$  we have
\begin{equation*}
\sum_{K \in{\cal{T}}_h} \kappa^m_K \sum_{\sigma' \in {\cal{E}}_K}
|\sigma'| \, \lambda^m_{\sigma'} \, (\V_h|_K \cdot \N_{K,\sigma'})
 =
\sum_{\sigma \in {\cal{E}}^{int}_h} |\sigma| \,
\phi_\sigma(\x_\sigma) \, (\Ut^m_h|_{L_\sigma}\cdot
\N_{K_\sigma,\sigma} - \Ut^m_h |_{K_\sigma} \cdot
\N_{K_\sigma,\sigma})=0,
\end{equation*}
it follows from (\ref{eq:mha}) that $(\Ut^m_h,\V)=(\F^m_h-\kappa^m_h
\nabla_h \widetilde p^m_h,\V)$ for any $\V \in \mathbf{RT_0}$. It
means that
\begin{equation*}
\Ut^m_h = \widetilde \Pi_{\mathbf{RT_0}} (\F^m_h - \kappa^m_h \,
\nabla_h \widetilde p^m_h) = \widetilde \Pi_{\mathbf{RT_0}} (\F^m_h
- \kappa^m_h \, \nabla_h p^m_h)=\U^m_h.
\end{equation*}
Thus $\U^m_h=\Ut^m_h$ satisfies (\ref{eq:mhc}) and $\hbox{div} \,
\U^m_h=s_h$.
 \qed


\section{Stability analysis}
\label{sec:stab}

\newcounter{sta}
\renewcommand{\theequation}{\thesection.\thesta}

We first check that a  maximum principle holds. It ensures that the
computed concentration and temperature are physically relevant.
\begin{prop}
\label{propmaxd} For any $m \in \{ 0, \dots,N \}$ we have
$0 \le c^m_{h}  \le 1$ and
 $\theta_- \le \theta^m_{h} \le \theta_+$.
\end{prop}
\noindent {\sc Proof.} We prove the result  by induction. Since
$c^0_h=\Pi_{P_0} c_0$ and $\theta^0_h=\Pi_{P_0} \theta_0$ the result
holds for $m=0$ thanks to  (\ref{eq:borncondi}) and
(\ref{eq:defprojp0}). Let us assume that  it is true for $m=n
\in\{0,\dots,N-1\}$. Let  $K \in {\cal{T}}_h$. Equation
(\ref{eq:eqcd})  implies
\begin{equation*}
\label{eq:inf1}
 (1+k \, s_K+ k  \, \lambda) \, c^{n+1}_K=c^n_K+k  \, s^{c}_K+k  \, D_c \sum_{\sigma \in {\cal{E}}_K \cap {\cal{E}}^{int}_h}
\tau_\sigma \,(c^{n+1}_{L_\sigma}-c^{n+1}_K) - k  \,
\bt_h(\U^n_h,c^{n+1}_h)\big|_K.
\end{equation*}
We consider the last term of this relation. Since for any $\sigma
\in {\cal{E}}_K \cap {\cal{E}}^{int}_h$  we have
\begin{equation*}
(\U^n_h \cdot \N_{K,\sigma})^+_\sigma \, c^{n+1}_K + (\U^n_h \cdot
\N_{K,\sigma})^-_\sigma
 \, c^{n+1}_{L_\sigma} =(\U^n_h \cdot \N_{K,\sigma})_\sigma \, c^{n+1}_K + (-\U^n_h \cdot
\N_{K,\sigma})^+_\sigma \, (c^{n+1}_K - c^{n+1}_{L_\sigma}).
\end{equation*}
We deduce from (\ref{eq:defbth})
$$- \bt_h(\U^n_h,c^{n+1}_h)\big|_K=
 \frac{1}{|K|} \Bigl(
-c^{n+1}_K \!\! \sum_{\sigma \in {\cal{E}}_K \cap {\cal{E}}^{int}_h}
\!\! |\sigma| \, (\U^n_h \cdot \N_{K,\sigma})_\sigma
+ \!\! \sum_{\sigma \in {\cal{E}}_K \cap {\cal{E}}^{int}_h} \!\!
(-\U^n_h \cdot \N_{K,\sigma})^+_\sigma
(c^{n+1}_{L_\sigma}-c^{n+1}_K) \Bigr).
$$
Since $\U^n_h \in  \mathbf{RT_0}$,  $(\U^n_h \cdot
\N_{K,\sigma})_\sigma=0$ for any $\sigma \in {\cal{E}}_K \cap
{\cal{E}}^{ext}_h$. It implies that $ \sum_{\sigma \in {\cal{E}}_K
\cap {\cal{E}}^{int}_h} |\sigma| \, (\U^n_h \cdot
\N_{K,\sigma})_\sigma =  \sum_{\sigma \in {\cal{E}}_K} |\sigma| \,
(\U^n_h \cdot \N_{K,\sigma})_\sigma$. Thus using the divergence
formula and proposition \ref{prop:propdivsh} we obtain
\begin{equation*}
\frac{1}{|K|} \sum_{\sigma \in {\cal{E}}_K \cap {\cal{E}}^{int}_h}
|\sigma| \, (\U^n_h \cdot \N_{K,\sigma})_\sigma =\frac{1}{|K|}
\int_K \hbox{\rm div}  \, \U^n_h \, d\x= s_K.
\end{equation*}
Therefore we get \stepcounter{sta}
\begin{align}
(1+2  \, k  \,s_K+ k  \, \lambda) \, c^{n+1}_K&=c^n_K+k  \,
s^{c}_K+k \, D_c \sum_{\sigma \in {\cal{E}}_K \cap
{\cal{E}}^{int}_h}
\tau_\sigma \, (c^{n+1}_{L_\sigma}-c^{n+1}_K) \nonumber \\
&+ \frac{1}{|K|}\sum_{\sigma \in {\cal{E}}_K \cap {\cal{E}}^{int}_h}
(-\U^n_h \cdot \N_{K,\sigma})^+_\sigma \,
(c^{n+1}_{L_\sigma}-c^{n+1}_K) . \label{eq:inf2}
\end{align}
We consider   $K_i \in {\cal{T}}_h$ such that $c^{n+1}_{K_i}=\min_{K
\in {\cal{T}}_h} c^{n+1}_K$. According to hypothesis
(\ref{eq:condpmax}) and definition (\ref{eq:defprojp0}) we have $2
\,  s_{K_i}+\lambda \ge 0$ and $s^c_{K_i} \ge 0$.
 Thus, using the induction hypothesis, we deduce from (\ref{eq:inf2})
\begin{equation*}
\min_{K \in {\cal{T}}_h} c^{n+1}_K=c^{n+1}_{K_i} \ge
\frac{c^n_{K_i}+k  s^{c}_{K_i}}{1+2  \, k \, s_{K_i}+k \,  \lambda}
\ge  \frac{k  \, s^{c}_{K_i}}{1+2\,  k  \, s_{K_i}+k \, \lambda} \ge
0.
\end{equation*}
We now consider $K_s \in {\cal{T}}_h$ such that
$c^{n+1}_{K_s}=\max_{K \in {\cal{T}}_h} c^{n+1}_K$. Using again
hypothesis (\ref{eq:condpmax})  we have $2 \,  s_{K_s}+\lambda \ge
s^c_{K_s} \ge 0$.
 Thus, using the induction hypothesis, we deduce from (\ref{eq:inf2})
 \begin{equation*}
\max_{K \in {\cal{T}}_h} c^{n+1}_K=c^{n+1}_{K_s} \le
\frac{c^n_{K_s}+k \,  s^{c}_{K_s}}{1+2  \, k \, s_{K_s}+k \,
\lambda} \le  \frac{1+k  \, s^{c}_{K_s}}{1+2\,  k  \, s_{K_s}+k \,
\lambda} \le 1.
\end{equation*}
A similar work for  equation (\ref{eq:eqtd})
 proves that $\theta_- \le \min_{K \in {\cal{T}}_h} \theta^{n+1}_K $ and $\max_{K \in {\cal{T}}_h} \theta^{n+1}_K \le \theta_+$.
Thus the induction hypothesis still holds for $m=n+1$. \qed

\noindent We now state the stability of the scheme
(\ref{eq:eqcd})-(\ref{eq:equd}).
\begin{prop}
\label{propstabs} For any $1 \le m \le N$ we have \stepcounter{sta}
\begin{eqnarray}
&& k \sum^m_{n=1} \|c^n_h\|^2_h +  k \sum^m_{n=1} \|\theta^n_h\|^2_h
\le C \, ,
\label{eq:estl2h1} \\
\stepcounter{sta} \label{eq:estud} && |\U^m_h| +|\nabla_h p^m_h|\le
C.
\end{eqnarray}
\end{prop}
\noindent {\sc Proof.} Let $0 \le n \le N-1$. Multiplying
(\ref{eq:eqcd}) by $2 \, k \, c^{n+1}_h$ we get
\begin{align*}
( c^{n+1}_h-c^n_h, 2 \, c^{n+1}_h )-2  \, k \, D_c\,  (\Delta_{h}
c^{n+1}_h,c^{n+1}_h)
  &+k \, \big((s_h+\lambda) \,
c^{n+1}_h,c^{n+1}_h\big) \nonumber \\
&+k  \,\b_h(\U^n_h,c^{n+1}_h,c^{n+1}_h) =k  \, (s_h^{c},c^{n+1}_h).
\end{align*}
We have $( c^{n+1}_h-c^n_h, 2 \,  c^{n+1}_h
)=|c^{n+1}_h|^2-|c^n_h|^2+|c^{n+1}_h-c^n_h|^2$. Thanks to
propositions \ref{prop:posbh} and \ref{prop:coerlap}
\begin{equation*}
-2 \, k \, (\Delta_{h} c^{n+1}_h,c^{n+1}_h)= 2 \,   k \,
\|c^{n+1}_h\|^2_h \, , \hspace{1cm}
 \b_h(\U^n_h,c^{n+1}_h,c^{n+1}_h) \ge 0.
\end{equation*}
Using  the Cauchy-Schwarz and Young inequalities we write
\begin{equation*}
k  \, (s_h^{c},c^{n+1}_h) \le k \, |s_h^{c}| \, |c^{n+1}_h| \le C \,
k \, |c^{n+1}_h|
 \le k \, \frac{\lambda}{2} \, |c^{n+1}_h|^2+C \, k.
\end{equation*}
Finally thanks to (\ref{eq:condpmax}) and (\ref{eq:defprojp0}) we
have $s_h \ge 0$. Thus we obtain
\begin{equation*}
|c^{n+1}_h|^2 - |c^n_h|^2+2 \, k \, D_c \, \|c^{n+1}_h\|^2_h +k \,
\frac{\lambda}{2} \, |c^{n+1}_h|^2 \le C \, k.
\end{equation*}
Let $m\in \{1,\dots,N\}$. Summing up the latter relation from $n=0$
to $m-1$  we get
\begin{equation*}
  |c^m_h|^2 + 2 \, k  \, D_c \sum_{n=1}^m \|c^n_h\|^2_h \le |c^0_h|^2
  +C  \sum_{n=1}^m k \le C \, ,
\end{equation*}
thanks to proposition \ref{propmaxd}. With a similar work on
equation (\ref{eq:eqtd}), we get (\ref{eq:estl2h1}). We now  prove
(\ref{eq:estud}). Let $n=m-1 \in \{0,\dots,N-1\}$. Multiplying
equation (\ref{eq:eqpd}) by $-p^{n+1}_h$ and using proposition
\ref{prop:propadjh}, we get \stepcounter{sta}
\begin{equation}
\label{eq:estpd5} (\kappa^{n+1}_h  \nabla_h p^{n+1}_h, \nabla_h
p^{n+1}_h)
  =(\F^{n+1}_h,\nabla_h p^{n+1}_h)+(\widetilde \Pi_{P^{nc}_1} s_h,p^{n+1}_h).
\end{equation}
The left-hand side term satisfies $  (\kappa^{n+1}_h  \nabla_h
p^{n+1}_h, \nabla_h p^{n+1}_h) \ge
  \kappa_{inf}  \, |\nabla_h p^{n+1}_h|^2$.
We now consider the right-hand side. Using  (\ref{eq:stabfh}), the
Cauchy-Schwarz and Young inequalities we write
 \begin{equation*}
\vert (\F^{n+1}_h,\nabla_h p^{n+1}_h) \vert \le |\F^{n+1}_h| \,
|\nabla_h p^{n+1}_h|
 \le  \frac{\kappa_{inf}}{4} \, |\nabla_h p^{n+1}_h|^2+C \, \|\F\|^2_{L^\infty(0,T;\L^2)}.
 \end{equation*}
 Also, the stability of $\widetilde \Pi_{P^{nc}_1}$ for the $L^2$-norm, proposition  \ref{proppoinp1nc} and the Young inequality lead to
 \begin{equation*}
\vert (\widetilde \Pi_{P^{nc}_1} s_h,p^{n+1}_h)\vert \le |s_h|\,
|p^{n+1}_h| \le C\, |p^{n+1}_h| \le C  \, |\nabla_h p^{n+1}_h| \le
\frac{\kappa_{inf}}{4}  \, |\nabla_h p^{n+1}_h|^2+ C.
 \end{equation*}
Thus we deduce from (\ref{eq:estpd5}) that
$ |\nabla_h p^{n+1}_h|^2=|\nabla_h p^{m}_h|^2 \le C$.
Then (\ref{eq:equd}) imply
\begin{equation*}
|\U^m_h|=|\U^{n+1}_h| \le |\F^{n+1}_h| + |\kappa^{n+1}_h \, \nabla_h
p^{n+1}_h| \le
\|\F\|_{L^\infty(0,T;\L^2)}+\|\kappa\|_{W^{1,\infty}((0,1) \times
(0,\infty))}\, |\nabla_h p^{n+1}_h| \le C.
\end{equation*}
Estimate (\ref{eq:estud}) is proven. \qed


\section{Convergence analysis}
\label{sec:conv}

\newcounter{cvg}
\renewcommand{\theequation}{\thesection.\thecvg}

Let $\eps=\max(h,k)$. In this section we study the behavior of the
scheme (\ref{eq:eqcd})-(\ref{eq:equd}) as $\eps \to 0$. We first
define the applications $c_\eps: \mathbb{R} \to P_0$, $\tilde
c_\eps: \mathbb{R} \to P_0$, $\theta_\eps: \mathbb{R} \to P_0$,
 $p_\eps: \mathbb{R} \to P^{nc}_1$, $s_\eps:\mathbb{R} \to P_0$, $s^c_\eps:\mathbb{R} \to P_0$ and  $\U_\eps: \mathbb{R} \to \mathbf{RT_0}$, $\F_\eps: \mathbb{R} \to \P_0$
 by setting
 for all $n\in\{ 0, \dots,N-1\}$ and $t\in [t_n,t_{n+1}]$
\begin{eqnarray*}
 && c_\eps(t)=c^{n+1}_h  , \hspace{0.5cm}  \tilde c_\eps(t)=c^n_h +
\frac{1}{k} \, (t-t_n) \, (c^{n+1}_h-c^n_h)  ,
 \hspace{0.5cm} \theta_\eps(t)=\theta^{n+1}_h  , \\
 && p_\eps(t)=p^{n+1}_h ,
  \hspace{0.5cm} s_\eps(t)=s_h , \hspace{.5cm} s^c_\eps(t)=s^c_h, \hspace{.5cm} \U_\eps(t)=\U_h^{n},
 \hspace{0.5cm} \F_\eps(t)=\F^{n+1}_h  ,
\end{eqnarray*}
and for all $t \not\in (0,T)$
\begin{equation*}
 c_\eps(t)=\tilde c_\eps(t)=\theta_\eps(t)=p_\eps(t)=s_\eps(t)=s_\eps^c(t)=0, \quad
\U_\eps(t)=\F_\eps(t)=\mathbf{0}.
\end{equation*}
We recall that the Fourier transform $\widehat{f}$ of a function $f
\in L^1(\mathbb{R})$ is defined for any $ \tau \in \mathbb{R}$ by
\stepcounter{cvg}
\begin{equation}
\label{eq:deftf}
 \widehat{f}(\tau)=\int_{\mathbb{R}} e^{-2i\pi \tau t} f(t) \, dt.
\end{equation}
We begin with the following estimate.
\begin{prop}
\label{prop:esttf} Let $0<\gamma<\frac{1}{4}$. There exists  $C>0$
such that for all $\eps>0$ 
\begin{equation*}
\label{eq:estcomptf}
 \int_{\mathbb{R}} |\tau|^{2  \gamma} \, (|\widehat c_\eps(\tau)|^2+\vert \widehat \theta_\eps(\tau)\vert^2) \, d\tau \le C.
\end{equation*}
\end{prop}
\noindent {\sc Proof.} Since  equations (\ref{eq:eqcd}) and
(\ref{eq:eqtd}) are similar we only prove the estimate on $\widehat
c_\eps$. We first define  $g_\eps:\mathbb{R} \to P_0 \cap L^2_0$ as
the solution of
\begin{equation*}
\label{eq:estcomp1}
 \Delta_{h} g_\eps= D_c \, \Delta_{h} c_\eps+s_\eps^{c} -(s_\eps+\lambda) \, c_\eps -  \bt_h(\U_\eps,c_\eps).
\end{equation*}
Multiplying this equation  by $-g_\eps$ we obtain \stepcounter{cvg}
\begin{equation}
\label{eq:esttf0}
 -\left(\Delta_{h} g_\eps,g_\eps\right)
 =-D_c  \big(\Delta_{h} c_\eps,g_\eps\big)
 -\big(s_\eps^{c} -(s_\eps+\lambda) \, c_\eps,g_\eps\big)
 +\b_h(\U_\eps,c_\eps,g_\eps).
\end{equation}
Proposition \ref{prop:coerlap} allows us to write 
\begin{equation*}
  -\left(\Delta_{h} g_\eps,g_\eps\right)=
 \|g_\eps\|^2_h, \hspace{2cm}  
 - \left(\Delta_{h} c_\eps,g_\eps\right) \le
 \|c_\eps\|_h \, \|g_\eps\|_h.
\end{equation*}
Thanks to the Cauchy-Schwarz inequality, (\ref{eq:stabfh})  and
proposition \ref{proppoinp0} we have
\begin{equation*}
\left|\big(s_\eps^{c} -(s_\eps+\lambda) \, c_\eps,g_\eps\big)\right|
 \le C \, (\vert s_c \vert + \vert s \vert + \lambda) \, |g_\eps|
  \le C \, \|g_\eps\|_h.
  \label{6.7}
\end{equation*}
According to proposition \ref{prop:stabbth}, then  proposition
\ref{propmaxd} and  (\ref{eq:estud}), we have
\begin{align*}
\left|\b_h\big(\U_\eps,c_\eps,g_\eps\big)\right|
 &\le C \,
 \|c_\eps\|_\infty \, \|g_\eps\|_h \,  \vert \mathrm{div} \, \U_\eps\vert
 +C \, \|c_\eps\|_h \, \|g_\eps\|_h \, \vert \U_\eps\vert \\
 &
 \le
 C \, \|g_\eps\|_h  \,  \vert \mathrm{div} \, \U_\eps\vert
 + C \, \|c_\eps\|_h \, \|g_\eps\|_h.
\end{align*}
Let us plug these estimates into  (\ref{eq:esttf0}) and integrate
from 0 to $T$. We get 
\begin{equation*}
\int^T_0 \|g_\eps\|_h \, dt \le C \int^T_0 \vert \mathrm{div} \,
\U_\eps \vert \, dt +C \int^T_0 \|c_\eps\|_h \, dt \le C \, ,
\end{equation*}
because of proposition \ref{prop:propdivsh} and (\ref{eq:estl2h1}).
Definition (\ref{eq:deftf}) then leads to \stepcounter{cvg}
\begin{equation}
\label{eq:estcomp5}
  \forall \, \tau \in \mathbb{R} \, , \hspace{1cm} \|\widehat g_\eps(\tau)\|_h \le C.
\end{equation}
We now use this estimate to prove (\ref{eq:estcomptf}). Equation
(\ref{eq:eqcd}) reads
\begin{equation*}
  \frac{d}{dt} \tilde c_\eps
  =\Delta_{h} g_\eps + (c^0_h \delta_0 - c^N_h  \delta_T)
\end{equation*}
where $\delta_0$ and $\delta_T$ are  Dirac distributions
respectively localized in $0$ and $T$. Let $\tau \in \mathbb{R}$.
Applying the Fourier transform to the latter equation we obtain
\begin{equation*}
  -2i\pi \tau \, \widehat{\tilde c}_\eps(\tau)
  =\Delta_{h} \widehat{g}_\eps(\tau) +(c^0_h - c^N_h e^{-2i\pi \tau T}).
\end{equation*}
Let us take the scalar product of this relation with $i
\hbox{sign}(\tau) \, \widehat{\tilde c}_\eps(\tau)$. Applying
propositions \ref{proppoinp0} and
\ref{prop:coerlap} leads to 
\begin{equation*}
\label{eq:estcomp16}
  2\pi  |\tau| \, |\widehat{\tilde c}_\eps(\tau)|^2
  \le C  \, \left( \|\widehat{g}_\eps(\tau)\|_h +|c^0_h|+|c^N_h|
  \right)
   \|\widehat{\tilde c}_\eps(\tau)\|_h.
\end{equation*}
We  assume that $\tau \neq 0$ and multiply this estimate by
$|\tau|^{2 \, \gamma-1}$. Using proposition \ref{propmaxd} and
(\ref{eq:estcomp5}) we get
\begin{equation*}
 |\tau|^{2 \, \gamma} \, |\widehat{\tilde c}_\eps(\tau)|^2 \le
 C \,
 |\tau|^{2 \, \gamma-1} \, \|\widehat{\tilde c}_\eps(\tau)\|_h.
\end{equation*}
Using the Young inequality and integrating over $\{\tau \in
\mathbb{R} \, ; |\tau| > 1\}$, we obtain
\begin{equation*}
\int_{|\tau| > 1} |\tau|^{2 \, \gamma} \, |\widehat{\tilde
c}_\eps(\tau)|^2 \, d\tau \le \int_{|\tau| > 1}  |\tau|^{4 \,
\gamma-2} \, d\tau + C \, \int_{|\tau| > 1} \|\widehat{\tilde
 c}_\eps(\tau)\|^2_h \, d\tau.
\end{equation*}
For $|\tau| \le 1$, we have $ |\tau|^{2 \, \gamma} \,
|\widehat{\tilde c}_\eps(\tau)|^2 \le |\widehat{\tilde
c}_\eps(\tau)|^2 \le C \, \|\widehat{\tilde c}_\eps(\tau)\|^2_h $
according to proposition \ref{proppoinp0}. Thus
\begin{equation*}
\int_{|\tau| \le  1} |\tau|^{2 \, \gamma} \, |\widehat{\tilde
c}_\eps(\tau)|^2 \, d\tau \le C \, \int_{|\tau| \le  1}
\|\widehat{\tilde c}_\eps(\tau)\|^2_h \, d\tau. 
\end{equation*}
By combining the bounds for $|\tau|>1$ and $|\tau| \le 1$ we get
\begin{equation*}
\int_{\mathbb{R}} |\tau|^{2 \, \gamma} \, |\widehat{\tilde
c}_\eps(\tau)|^2 \, d\tau \le \int_{|\tau| > 1}  |\tau|^{4 \,
\gamma-2} \, d\tau+ C \, \int_{\mathbb{R}} \|\widehat{\tilde
c}_\eps(\tau)\|^2_h \, d\tau.
\end{equation*}
Since $4 \, \gamma-2<-1$, we have $\int_{|\tau| > 1}  |\tau|^{4 \,
\gamma-2} \, d\tau \le C$. Thanks to the Parseval theorem and
(\ref{eq:estl2h1})
\begin{equation*}
  \int_{\mathbb{R}} \|\widehat{\tilde
 c}_\eps(\tau)\|^2_h \, d\tau
\le  \int_{\mathbb{R}} \|{\tilde
 c}_\eps\|^2_h \, dt
 \le C \, \left( k \, \|c^0_h\|^2_h+k \sum_{n=1}^N \|c^n_h\|^2_h \right)\le C \, ,
\end{equation*}
because $\|c^0_h\|_h=\|\Pi_{P_0} c_0\|_h \le C \, \|c_0\|_1$ (see
\cite{eym} p. 776). Hence the result.  \qed

\vspace{.2cm}

\noindent We can now prove the following convergence result.
\begin{prop}
There exists a subsequence of
$(c_\eps,\theta_\eps,p_\eps,\U_\eps)_{\eps>0}$, not relabeled for
convenience, such that the following convergences hold  for  $\eps
\to 0$ \stepcounter{cvg}
 \begin{eqnarray}
 \label{eq:propconv1}
&&   c_\eps \to c \  \hbox{ in } L^2(0,T;L^2)  , \hspace{2cm}
   \theta_\eps \to \theta \   \hbox{ in } L^2(0,T;L^2)   ,
\\
\stepcounter{cvg} &&   p_\eps \rightharpoonup p \  \hbox{ weakly in
} L^2(0,T;H^1) , \hspace{.75cm}
   \U_\eps \rightharpoonup \U \  \hbox{ weakly in } L^2(0,T;\L^2).
\label{eq:propconv2}
 \end{eqnarray}
The limits $(c,\theta,p,\U)$ satisfy the following properties.  We
have
  $c \in L^2(0,T;H^1)$, $\theta \in L^2(0,T;H^1)$,
  $p \in L^\infty(0,T;H^1)$ and
$\U \in L^\infty(0,T;\L^2)$. We also have $0\le c(\x,t) \le 1$ and
$\theta^- \le \theta(\x,t) \le \theta^+$ a.e. in $\Omega \times
[0,T]$.
For all $\phi \in {\mathcal{C}}^\infty_0(\Omega \times (-1,T))$, $c$
and  $\theta$ satisfy \stepcounter{cvg}
\begin{eqnarray}
\label{eq:propconv3}
 && \int_0^T \Big( (c,\partial_t \phi)+D_c \, (\nabla c,\nabla \phi) - c\, (\U \cdot \nabla \phi)
  -\big(s_c-(s+\lambda) \, c \big) \, \phi \Big) \,
  dt=\big(c_0,\phi(\cdot,0)\big) \, , \\
\stepcounter{cvg} \label{eq:propconv4} &&  \int_0^T \Big(
(\theta,\partial_t \phi)+D_\theta \, (\nabla \theta,\nabla \phi) -
\theta\, (\U \cdot \nabla \phi)
  +\big(s_\theta+s \, (\theta-\theta_*) \big) \, \phi \Big) \,
  dt=\big(\theta_0,\phi(\cdot,0)\big).
\end{eqnarray}
Lastly we have \stepcounter{cvg}
\begin{equation}
 \U=\F-\kappa(c,\theta) \, \nabla p  \quad \hbox{ in }
L^\infty(0,T;\L^2) \, , \hspace{2cm}  \hbox{\rm div} \, \U=s \quad
\hbox{ in } L^2.
\end{equation}
\end{prop}
\noindent {\sc Proof.} In what follows, the convergence results hold
for extracted subsequences. They are not relabeled   for
convenience. We begin by proving (\ref{eq:propconv1}). According to
proposition \ref{propmaxd}, the sequence $(c_\eps)_{\eps>0}$ is
uniformly bounded in $L^\infty(0,T;L^2)$. Thus there exists  $c \in
L^\infty(0,T;L^2)$ such that
\begin{equation*}
\label{eq:convfor1}
 c_\eps \rightharpoonup c \  \hbox{ weakly in } \, L^2(0,T;L^2).
\end{equation*}
Using the Fourier transform, we  prove that this convergence is
strong. Let $d_\eps=c_\eps-c$ and $M>0$. We use the following
splitting \stepcounter{cvg}
\begin{equation}
\label{eq:convfor2}
  \int_{\mathbb{R}} |\widehat{d}_\eps(\tau)|^2 \, d\tau
  =
   \int_{|\tau| > M} |\widehat d_\eps(\tau)|^2 \, d\tau
 + \int_{|\tau| \le M} |\widehat d_\eps(\tau)|^2 \, d\tau
  =I^M_\eps+J^M_\eps.
\end{equation}
Since $ |\widehat d_\eps(\tau)|^2 \le 2 |\widehat c_\eps(\tau)|^2+2
|\widehat c(\tau)|^2$ we have
\begin{equation*}
  I^M_\eps
  \le 2 \int_{|\tau| > M} |\widehat c_\eps(\tau)|^2 \, d\tau +
2 \int_{|\tau| > M} |\widehat c(\tau)|^2 \, d\tau.
\end{equation*}
Using proposition \ref{prop:esttf} we write
\begin{equation*}
  \int_{|\tau| > M} |\widehat c_\eps(\tau)|^2 \, d\tau
  \le \frac{1}{M^{2  \gamma}}  \int_{|\tau| > M} |\tau|^{2  \gamma} \, |\widehat c_\eps(\tau)|^2 \, d\tau
 \le \frac{C}{M^{2 \gamma}}.
\end{equation*}
Hence
\begin{equation*}
I^M_\eps \le \frac{2  C}{M^{2  \gamma}} +2 \int_{|\tau| > M}
|\widehat c(\tau)|^2 \, d\tau.
\end{equation*}
 This implies that for all $\eps>0$, $I^M_\eps \to 0$ when $M \to \infty$.
We now consider $J^M_\eps$. Let $\tau \in [-M,M]$.  Since $c_\eps(t)
\in P_0$ for all $t \in \mathbb{R}$, and  $c_\eps \rightharpoonup c$
weakly in $L^2(0,T;L^2)$, we deduce from (\ref{eq:deftf}) that
 $\widehat c_\eps(\tau) \in P_0$ and $\widehat c_\eps(\tau) \rightharpoonup \widehat c(\tau)$ weakly in $L^2$.
Extanding $\widehat c_\eps(\tau)$ by $0$ outside $\Omega$, one
checks (\cite{eym} p.811) that
\begin{equation}
\label{eq:convfor201}
 \forall \,
\mathbf{\eta} \in \mathbb{R}^2 \, , \hspace{1cm} |\widehat
c_\eps(\tau)(\cdot+\mathbf{\eta})-\widehat c_\eps(\tau)| \le C \,
 \|\widehat c_\eps(\tau)\|_h \, |\mathbf{\eta}| \, (|\mathbf{\eta}|+h).
\end{equation}
Then, using estimate (\ref{eq:estl2h1}), we deduce from \cite{eym}
(p.834) that $\widehat c_\eps(\tau) \to \widehat c(\tau)$ strongly
in $L^2$. Thus $\widehat d_\eps(\tau)=\widehat c_\eps(\tau)-\widehat
c(\tau) \to 0$ in $L^2$, so that $J^\eps_M \to 0$ when $\eps \to 0$.
Now, let us report the limits for $I^\eps_M$ and $J^\eps_M$ into
(\ref{eq:convfor2}). Using the Parseval identity we get
\begin{equation*}
  \int_{\mathbb{R}} |\widehat d_\eps(\tau)|^2 \, d\tau
  =\int_{\mathbb{R}} |d_\eps|^2 \, dt
  =\int_{\mathbb{R}} |c_\eps-c|^2 \, dt \to 0.
\end{equation*}
Thus we have proven that $c_\eps \to c$ in $L^2(0,T;L^2)$.
 A similar work
proves  that $ \theta_\eps \to \theta$  in $L^2(0,T;L^2)$ with
$\theta \in L^\infty(0,T;L^2)$. Hence (\ref{eq:propconv1}) is
proven. Moreover using proposition \ref{propmaxd} we obtain  $0 \le
c(\x,t) \le 1$ and $\theta^- \le \theta(\x,t) \le \theta^+$ a.e. in
$\Omega \times [0,T]$. Lastly, using
 (\ref{eq:estl2h1}) and (\ref{eq:convfor201}), we get   as in \cite{eym} (p.811)
 that $c \in L^2(0,T;H^1)$ and $\theta \in L^2(0,T;H^1)$.

\noindent Let us now consider the sequences $(p_\eps)_{\eps>0}$ and
$(\U_\eps)_{\eps>0}$.  According to (\ref{eq:defprojp1c}) and
(\ref{eq:estud}) the sequence $(\Pi_{P^c_1} p_\eps)_{\eps>0}$ is
bounded in $L^\infty(0,T;H^1)$. It implies that there  exists $p \in
L^\infty(0,T;H^1)$ such that $\Pi_{P^c_1}p_\eps \rightharpoonup p$
weakly in $L^2(0,T;H^1)$. Using proposition \ref{prop:errint} we get
$p_\eps \rightharpoonup p$ weakly in $L^2(0,T;H^1)$. Moreover,
according to (\ref{eq:estud}),
  the sequence $(\U_\eps)_{\eps>0}$ is
bounded in $L^\infty(0,T;\L^2)$ . Thus we  have $\U_\eps
\rightharpoonup \U$ weakly in $L^2(0,T;\L^2)$ with $\U \in
L^\infty(0,T;\L^2)$. We check the properties of $\U$.
 Using a Taylor expansion, the  Cauchy-Schwarz inequality, and  a
 density argument, we have
\begin{equation*}
  \|\kappa(c,\theta)-\kappa(c_\eps,\theta_\eps)\|_{L^2(0,T;L^2)}
  \le \|\kappa\|_{W^{1,\infty}((0,1) \times (0,\infty))} \, (\|c-c_\eps\|_{L^2(0,T;L^2)}
 + \|\theta-\theta_\eps\|_{L^2(0,T;L^2)}).
\end{equation*}
Thus, using the strong convergence of the sequences
$(c_\eps)_{\eps>0}$ and $(\theta_\eps)_{\eps>0}$, we have
  $\kappa(c_\eps,\theta_\eps) \to \kappa(c,\theta)$ in $L^2(0,T;L^2)$.
Since $\nabla_h p_\eps \rightharpoonup \nabla p$ weakly in
$L^2(0,T;\L^2)$, we deduce from this \stepcounter{cvg}
\begin{equation}
\label{eq:convfor26} \kappa(c_\eps,\theta_\eps) \,   \nabla_h p_\eps
\rightharpoonup
 \kappa(c,\theta)  \nabla p \hspace{.5cm} \hbox{ weakly in }L^2(0,T;\L^2).
\end{equation}
Now let $\V \in L^2(0,T;({\mathcal{C}}^\infty_0)^2)$. According to
(\ref{eq:equd}) we have
\begin{equation*}
  (\U_\eps,\Pi_{\mathbf{RT_0}}\V)=(\F_\eps-\kappa(c_\eps,\theta_\eps) \, \nabla_h p_\eps,\V).
\end{equation*}
Using proposition \ref{prop:errint}  one checks easily that
$(\F_\eps,\Pi_{\mathbf{RT_0}}\V) \to (\F,\V)$
 and $(\U_\eps,\Pi_{\mathbf{RT_0}}\V) \to (\U,\V)$  in
$L^1(0,T)$. Using moreover convergence (\ref{eq:convfor26}) and a
density argument, we deduce from this that $\U=\F-\kappa(c,\theta)
\, \nabla p$. And since $\hbox{div}  \, \U_\eps=s_\eps$ by
proposition \ref{prop:propdivsh}, we also have $\hbox{div} \, \U=s$.

\noindent We finally prove that $c$ satisfies (\ref{eq:propconv3}). 
For all $t \in (0,T)$ equation (\ref{eq:eqcd}) reads
\begin{equation*}
 \frac{d}{dt}\tilde c_\eps - D_c \, \Delta_{h}  c_\eps +\bt_h(\U_\eps,c_\eps)=s_\eps^{c} - (s_\eps+\lambda) \, c_\eps.
\end{equation*}
Let  $\psi \in {\cal{C}}^{\infty}_0(\Omega \times (-1,T))$ and
$\psi_h=\widetilde \Pi_{P_0} \psi$. Multiplying the latter equation
by $\psi_h$ and integrating over $[0,T]$, we obtain
\stepcounter{cvg}
\begin{equation}
 \int^T_0  \! \left(\frac{d}{dt}\tilde
c_\eps,\psi_h\right) \, dt
  - D_c \int^T_0 \! (\Delta_{h} c_\eps,\psi_h) \, dt
  + \int^T_0 \! \b_h(\U_\eps,c_\eps,\psi_h) \, dt
= \int^T_0 \left(s_\eps^{c}-(s_\eps+\lambda) \, c_\eps,\psi_h\right)
\, dt. \label{eq:formvar35}
\end{equation}
We now pass to the limit $\eps \to 0$ in this equation. We begin
with the term $\int^T_0 \b_h(\U_\eps,c_\eps,\psi_h)  \, dt$. We use
the splitting
$ \b(\U,c,\psi)- \b_h(\U_\eps,c_\eps, \psi_h) =A^\eps_1+
A^\eps_2+A^\eps_3$
with
\begin{eqnarray*}
&& A^\eps_1 = \b(\U,c,\psi)-\b(\U_\eps,c,\psi)  , \qquad A^\eps_2=
\b(\U_\eps,c,\psi)-\int_\Omega \hbox{div}(c  \, \U_\eps) \, \psi_h
\, d\x \ ,
\\
&& A^\eps_3= \int_\Omega \hbox{div}(c  \, \U_\eps) \, \psi_h \, d\x
- \b_h(\U_\eps,c_\eps,\psi_h).
\end{eqnarray*}
According to definition (\ref{eq:defb})
\begin{equation*}
  A^\eps_1=\b(\U,c,\psi)-\b(\U_\eps,c,\psi)
  =-\int_\Omega c \, (\U-\U_\eps) \cdot  \nabla \psi \, d\x.
\end{equation*}
We know  that $c  \, \nabla \psi \in L^2(0,T;\L^2)$. Since $\U_\eps
\rightharpoonup \U$ in $L^2(0,T;\L^2)$ we get
$ \int^T_0 A^\eps_1 \, dt \to 0$.
We now consider $A^\eps_2$. We  have 
\begin{equation*}
  A^\eps_2=\int_\Omega (\psi-\psi_h) \, \hbox{div}(c \, \U_\eps) \, d\x
  =\int_\Omega (\psi-\psi_h) \, (\U_\eps \cdot \nabla c+ c \, \hbox{div}  \,  \U_\eps) \, d\x.
\end{equation*}
Using the Cauchy-Schwarz inequality we get
\begin{equation*}
\label{eq:formvar405}
  \int^T_0 |A^\eps_2| \, dt \le \|\psi-\psi_h\|_{L^\infty(\Omega \times (0,T))} \,
  (\|\U_\eps\|_{L^2(0,T;\L^2)}+\|\hbox{div} \, \U_\eps\|_{L^2(0,T;\L^2)}) \, \|c\|_{L^2(0,T;H^1)}.
\end{equation*}
Using a Taylor expansion, one checks that
$   \|\psi-\psi_h\|_{L^\infty(\Omega \times (0,T))}
   \le  h \, \|\nabla \psi\|_{L^\infty(\Omega \times (0,T))} $.
Thus
$ \int^T_0 A^\eps_2 \, dt \to 0$. Finally we estimate $A^\eps_3$.
For all triangles $K \in {\mathcal{T}}_h$ and $L\in {\mathcal{T}}_h$
sharing an edge $\sigma$, we set $c_{K,L}=c_K$ if $\U_\eps \cdot
\N_{K,\sigma} \ge 0$ and $c_{K,L}=c_L$ otherwise. Using the
divergence formula, we deduce from  definition (\ref{eq:defbh})
\begin{eqnarray*}
A^\eps_3 &=& \sum_{K \in {\cal{T}}_h} \sum_{\sigma \in {\cal{E}}_K
\cap {\cal{E}}^{int}_h}
  \psi_K \int_\sigma (c-c_{K,L_\sigma}) \, (\U_\eps \cdot \N_{K,\sigma}) \, d\sigma
\\
&=&  \sum_{\sigma \in {\cal{E}}^{int}_h}
  (\psi_{K_\sigma}-\psi_{L_\sigma}) \int_\sigma (c-c_{K_\sigma,L_\sigma}) \, (\U_\eps \cdot \N_{K_\sigma,\sigma})\, d\sigma.
\end{eqnarray*}
Using definition (\ref{eq:defprojp1nc}) this reads
\begin{eqnarray*}
  A^\eps_3&=& \sum_{\sigma \in {\cal{E}}^{int}_h}
  (\psi_{K_\sigma}-\psi_{L_\sigma}) \int_\sigma (\Pi_{P^{nc}_1}c-c_{K_\sigma,L_\sigma}) \, (\U_\eps \cdot \N_{K_\sigma,\sigma})\, d\sigma \\
&=& \sum_{\sigma \in {\cal{E}}^{int}_h}
  (\psi_{K_\sigma}-\psi_{L_\sigma}) \,  |\sigma| \left( (\Pi_{P^{nc}_1}c)(\x_\sigma)-c_{K_\sigma,L_\sigma}\right)  (\U_\eps \cdot \N_{K_\sigma,\sigma})_\sigma.
\end{eqnarray*}
Using a  Taylor expansion, one checks that
$|\psi_{K_\sigma}-\psi_{L_\sigma}| \le h \, \|\nabla
\psi\|_{L^\infty(\Omega \times (0,T))}$. Moreover $|\sigma| \le h$.
Thus, using the Cauchy-Schwarz inequality, we have
\begin{eqnarray*}
  |A^\eps_3| &\le&
   C  \, h^2 \sum_{\sigma \in {\cal{E}}^{int}_h} |\U_\eps(\x_\sigma)|
   \left|(\Pi_{P^{nc}_1} c)(\x_\sigma) -c_{K_\sigma,L_\sigma}\right|
\\
&\le&
 C  \, h^2
 \bigl(\sum_{\sigma \in {\cal{E}}^{int}_h} |\U_\eps(\x_\sigma)|^2 \bigr)^{1/2}
\bigl(\sum_{\sigma \in {\cal{E}}^{int}_h}  |(\Pi_{P^{nc}_1}
c)(\x_\sigma) -c_{K_\sigma,L_\sigma}|^2 \bigr)^{1/2}.
\end{eqnarray*}
Using the assumption on the mesh, one checks that $|K| \ge C \, h^2$
for all $K \in {\mathcal{T}}_h$. Thus, thanks to a   quadrature
formula, we have
\begin{eqnarray*}
  \hspace*{-1.1cm}
 |A^\eps_3|
&\le& C \Bigl( \sum_{K \in {\cal{T}}_h} \frac{|K|}{3} \!\!\!\!
\sum_{\sigma \in {\cal{E}}_K \cap {\cal{E}}^{int}_h} \!\!
|\U_\eps(\x_\sigma)|^2 \Bigr)^{1/2} \Bigl( \sum_{K \in {\cal{T}}_h}
\frac{|K|}{3} \!\!\!\!  \sum_{\sigma \in {\cal{E}}_K \cap
{\cal{E}}^{int}_h} \!\! |(\Pi_{P^{nc}_1} c)(\x_\sigma)
-c_{K_\sigma}|^2 \Bigr)^{1/2}
\\
&\le& C \, |\U_\eps| \, |\Pi_{P^{nc}_1} c -c_\eps|.
\end{eqnarray*}
We write $\Pi_{P^{nc}_1} c -c_\eps= ( \Pi_{P^{nc}_1} c -c ) +
(c-c_\eps)$ and  we use proposition \ref{prop:errint} . We obtain
with (\ref{eq:estud})
\begin{equation*}
\int^T_0 |A^\eps_3| \, dt \le C \, \|\U_\eps\|_{L^\infty(0,T;\L^2)}
\, ( h \, \|c\|_{L^2(0,T;H^1)}+\|c-c_\eps\|_{L^2(0,T;L^2)}).
\end{equation*}
Since $c_\eps \to c$  in $L^2(0,T;L^2)$ when $\eps=\max(h,k) \to 0$,
we conclude that
 $ \int^T_0 A^\eps_3 \, dt \to 0$.
Gathering the limits for $A^\eps_1$, $A^\eps_2$, $A^\eps_3$, we
obtain
\begin{equation*}
\label{eq:formvar435}
  \int^T_0 \b_h(\U_\eps,c_\eps, \widetilde \Pi_{P_0} \psi) \, dt \to \int^T_0 \b(\U,c,\psi) \, dt.
\end{equation*}
We now consider the other terms in (\ref{eq:formvar35}). Proposition
\ref{prop:adjlap} leads to \stepcounter{cvg}
\begin{eqnarray}
&&  (\Delta_{h} c_\eps,\widetilde \Pi_{P_0} \psi)
  =\big( c_\eps, \Delta_{h} ( \widetilde \Pi_{P_0} \psi) \big)
  =\big( c_\eps, \Delta_{h} ( \widetilde \Pi_{P_0} \psi) -\Delta \psi \big)+
  (c_\eps,\Delta \psi).
\label{eq:formvar44}
\end{eqnarray}
According to proposition \ref{propconslap}
\begin{equation*}
\left|\big( c_\eps, \Delta_{h} ( \widetilde \Pi_{P_0} \psi) -\Delta
\psi\big)\right|
 \le \|c_\eps\|_h \, \|\Delta_{h} ( \widetilde \Pi_{P_0} \psi) - \Delta \psi \|_{-1,h}
 \le C \, h \, \|c_\eps\|_{h} \, \|\psi\|_2.
\end{equation*}
We then apply the  Cauchy-Schwarz inequality and use
(\ref{eq:estl2h1}). We obtain
\begin{equation*}
\int^T_0 \left|\big( c_\eps, \Delta_{h} ( \widetilde \Pi_{P_0} \psi)
-\Delta \psi \big)\right| \, dt \le C  \, h \, \Bigl( \int^T_0
\|c_\eps\|^2_h \, dt\Bigr)^{1/2} \le C  \, h \, \Bigl( k
\sum_{n=1}^N \|c^n_h\|^2_h \Bigr)^{1/2} \le C  \, h.
\end{equation*}
Moreover, since $c_\eps \to c$ in $L^2(0,T;L^2)$, we have $\int^T_0
(c_\eps,\Delta \psi)\,   dt
 \to \int^T_0 (c,\Delta \psi) \,  dt$.
Thus we deduce from (\ref{eq:formvar44})
\begin{equation*}
\label{eq:formvar45} \int^T_0 (\Delta_{h} c_\eps,\widetilde
\Pi_{P_0} \psi) \, dt \to \int^T_0 (c,\Delta \psi) \, dt.
\end{equation*}
We are left with two terms. First, using Taylor expansions, one
checks that \stepcounter{cvg}
\begin{equation}
\label{eq:convpsih} \psi_h \to \psi \, , \hspace{.2cm}
 \partial_t \psi_h \to \partial_t \psi \hspace{.2cm} \hbox{ in } L^2(\Omega \times (-1,T))
  \, , \hspace{1cm} \psi_h(\cdot,0) \to \psi(\cdot,0) \hbox{ in }L^2.
\end{equation}
We know that $c_\eps \to c$ in $L^2(0,T;L^2)$. Thus
\begin{equation*}
\label{eq:formvar46}
 \int^T_0 \left( s_{c}-(s+\lambda) \, c_\eps,\psi_h\right) \, dt
 \to   \int^T_0 \left( s_c-(s+\lambda) \, c,\psi\right) \, dt.
\end{equation*}
Finally, integrating by parts the first term of
(\ref{eq:formvar35}), we get
\begin{equation*}
  \int^T_0 \Bigl(\frac{d}{dt} \tilde c_\eps, \psi_h \Bigr) \, dt
  =(\tilde c_\eps, \psi_h)_{t=T} -(\tilde c_\eps, \psi_h)_{t=0}
  -\int^T_0 (\tilde c_\eps, \partial_t \psi_h) \, dt.
\end{equation*}
Since $\psi \in {\cal{C}}^\infty_0(\Omega \times (-1,,T))$ we have
$(\tilde c_\eps, \psi_h)_{t=T}=0$. Using proposition
\ref{prop:errint} one checks that $c^0_h=\Pi_{P_0}c_0 \to c_0$ in
$L^2$; using moreover (\ref{eq:convpsih})  we get
\begin{equation*}
(\widetilde c_\eps,  \psi_h)_{t=0}= \big(c^0_h, \psi_h(\cdot,0)\big)
=\big(\Pi_{P_0} c_0,  \psi_h(\cdot,0)\big) \to
\big(c_0,\psi(\cdot,0)\big).
\end{equation*}
For the last term, one easily checks  that $\|\tilde
c_\eps-c_\eps\|_{L^2(0,T;L^2)} \to 0$. Thus, since $c_\eps \to c$ in
$L^2(0,T;L^2)$, we also have  $\tilde c_\eps \to c$ in
$L^2(0,T;L^2)$. Using moreover (\ref{eq:convpsih}) we get
 $\int^T_0 (\tilde c_\eps,  \partial_t \psi_h) \,   dt
 \to \int^T_0  (c  , \partial_t \psi) \, dt$.
 Therefore
\begin{equation*}
\label{eq:formvar47}
 \int^T_0 \Bigl(\frac{d}{dt} \tilde c_\eps, \psi_h \Bigr) \, dt
 \to -\big(c_0,\psi(\cdot,0)\big)-\int^T_0 (c, \partial_t \psi) \, dt.
\end{equation*}
By gathering the limits we have obtained in (\ref{eq:formvar35}) we
get (\ref{eq:propconv3}). The relation (\ref{eq:propconv4})  for
$\theta$ is proven in a similar way.
 \qed


\section{Error estimates}
\label{sec:esterr}

\newcounter{err}
\renewcommand{\theequation}{\thesection.\theerr}

 We have proven in section \ref{sec:conv} that the problem (\ref{eq:1.10})--(\ref{eq:1.14}) has a weak solution
  $(c,\theta,p,\U)$. From now on, we assume the following regularity for this solution:
\begin{eqnarray*}
  && c,\theta \in {\cal{C}}(0,T;H^2)   \, ,
  \hspace{1cm} c_t,\theta_t \in L^2(0,T;H^{1+r})\cap {\cal{C}}(0,T;L^2) \, , \\
&&   c_{tt},\theta_{tt} \in L^2(0,T;L^2) \, , \hspace{.5cm}
 p \in {\cal{C}}(0,T;H^2)  \, , \hspace{.5cm} \U \in {\cal{C}}(0,T;\H^{1+s}) \, ,
\end{eqnarray*}
with $r>0$ and $s>0$. We also assume that $\F \in
{\mathcal{C}}(0,T;\H^1)$.

\subsection{Definitions}

For all $m\in\{0,\dots,N\}$, we define the following errors
  \begin{align*}
    e^m_{h,c}&= c(t_m) -c^m_h  ,
    \qquad 
     e^m_{h,\theta}= c(t_m) -\theta^m_h  ,
     \\
    e^m_{h,p}&= p(t_m) -p^m_h  , \qquad     
    \mathbf{e}^m_{h,\U}= \U(t_m) -\U^m_h.
  \end{align*}
We have the following splittings
\begin{eqnarray*}
  && e^m_{h,c}=\eps^m_{h,c}+\eta^m_{h,c}  ,
  \qquad
  e^m_{h,\theta}=\eps^m_{h,\theta}+\eta^m_{h,\theta}  , \\
  && e^m_{h,p}=\eps^m_{h,p}+\eta^m_{h,p}  ,
  \qquad
  \mathbf{e}^m_{h,\U}=\Eps{m}{h,\U}+\Eta{m}{h,\U} ,
\end{eqnarray*}
with the discrete errors
  \begin{align*}
    \eps^m_{h,c} &=  \widetilde \Pi_{P_0} c(t_m)-c^m_h ,
    \qquad
       \eps^m_{h,\theta} =  \widetilde \Pi_{P_0} \theta(t_m)-\theta^m_h   ,   \\ 
       \eps^m_{h,p} &=  \Pi_{P^{nc}_1} p(t_m)- p^m_h  ,
    \quad \, \,
    \Eps{m}{h,\U} =   \Pi_{\mathbf{RT_0}} \U(t_m)-\U^m_h , 
  \end{align*}
and the interpolation errors
  \begin{align*}
     \eta^m_{h,c} &= c(t_m)-\widetilde \Pi_{P_0} c(t_m) ,
     \qquad
      \eta^m_{h,\theta} = \theta(t_m)-\widetilde \Pi_{P_0} \theta(t_m)
      ,      \\
    \eta^m_{h,p} &= p(t_m)-\Pi_{P^{nc}_1} p(t_m)  , \quad \, \, 
     \Eta{m}{h,\U} =  \U(t_m)-\Pi_{\mathbf{RT_0}} \U(t_m). 
  \end{align*}
The interpolation errors are estimated as follows. We write
  $|\eta^m_{h,c}| \le  |c(t_m)-\Pi_{P_0} c(t_m)|+|\Pi_{P_0} c(t_m)-\widetilde \Pi_{P_0}
  c(t_m)|$
and the same for $\eta^m_{h,\theta}$. Using proposition
\ref{prop:errint} and (\ref{eq:esti}) we obtain \stepcounter{err}
\begin{equation}
\label{eq:esterrintc}
  |\eta^m_{h,c}| \le C \,   h  \, \|c(t_m)\|_1 \le C  \, h  \, \|c\|_{L^\infty(0,T;H^1)} \, ,
  \hspace{1cm}
   |\eta^m_{h,\theta}|  \le C  \, h  \,\|\theta\|_{L^\infty(0,T;H^1)}.
\end{equation}
According to  proposition \ref{prop:errint} we also have
\stepcounter{err}
\begin{eqnarray}
\label{eq:esterrintp} && |\eta^m_{h,p}|+|\widetilde \nabla_h
\eta^m_{h,p}| \le C  \, h \, \|p(t_m)\|_2 \le C  \, h \,
\|p\|_{L^\infty(0,T;H^2)} ,
\\
\stepcounter{err} \label{eq:esterrintu} && |\eta^m_{h,\U}| \le C \,
h \, \|\U(t_m)\|_1 \le C \, h  \, \|\U\|_{L^\infty(0,T;\H^1)}.
\end{eqnarray}
We now have to estimate the discrete errors.

\begin{prop}
\label{propeqerr} For all $n\in\{0,\dots,N-1\}$ and  $\psi_h \in
P^{nc}_1$ we have \stepcounter{err}
\begin{eqnarray}
\hspace{-1.2cm}&& \frac{\eps^{n+1}_{h,c} -\eps^n_{h,c}}{k} - D_c \,
\Delta_h \eps^{n+1}_{h,c} + \bt_h\big(\Eps{n}{h,\U},\widetilde
\Pi_{P_0} c(t_{n+1})\big) + \bt_h(\U^n_h,\eps^{n+1}_{h,c})
+(s_h^{n}+\lambda) \, \eps^{n+1}_{h,c}= C^{n+1}_{h,1}+ C^{n+1}_{h,2}
,
\label{eq:eqerr26} \\
\stepcounter{err} \hspace{-1.2cm}&&  \frac{\eps^{n+1}_{h,\theta}
-\eps^n_{h,\theta}}{k} - D_\theta \, \Delta_{h}
\eps^{n+1}_{h,\theta} + \bt_h\big(\Eps{n}{h,\U},\widetilde \Pi_{P_0}
\theta(t_{n+1})\big) + \bt_h(\U^n_h,\eps^{n+1}_{h,\theta})
+ s_h \, \eps^{n+1}_{h,\theta}=
\Theta^{n+1}_{h,1}+\Theta^{n+1}_{h,2} ,
\label{eq:eqerr27} \\
\stepcounter{err} \hspace{-1.2cm}&&
\big(\kappa(c^{n+1}_h,\theta^{n+1}_h) \, \nabla_h
\eps^{n+1}_{h,p},\nabla_h \psi_h\big) = -\big( ( \kappa^{n+1}_{h,1}
\eps^{n+1}_{h,c}+
  \kappa^{n+1}_{h,2} \eps^{n+1}_{h,\theta})  \, \nabla p(t_{n+1}) ,\nabla_h \psi_h \big)
- \big(\P^{n+1}_{h},\nabla_h \psi_h \big) ,
 \label{eq:eqerr28} \\
\stepcounter{err} \hspace{-1.2cm}&& \Eps{n+1}{h,\U} = -\widetilde
\Pi_{\mathbf{RT_0}} \big( (\kappa^{n+1}_{h,1}  \eps^{n+1}_{h,c}+
  \kappa^{n+1}_{h,2}  \eps^{n+1}_{h,\theta} ) \, \nabla p(t_{n+1})
  +\kappa(c^{n+1}_h,\theta^{n+1}_h) \,
   \nabla_h \eps^{n+1}_{h,p}\big)
- \mathbf{U}^{n+1}_h. \label{eq:eqerr29}
\end{eqnarray}
For all $m\in\{0,\dots,N\}$, the consistency errors $C^m_{h,1}$,
$C^m_{h,2}$, $\Theta^m_{h,1}$, $\Theta^m_{h,2}$, $\P^m_h$ and
$\mathbf{U}^m_h$ are defined in {\rm (\ref{eq:deferrconsc1})}, {\rm
(\ref{eq:deferrconstheta1})}, {\rm (\ref{eq:deferrconsp})},   {\rm
(\ref{eq:deferrconsu})} and the terms $\kappa^m_{h,1}$ and
$\kappa^m_{h,2}$ are given by {\rm (\ref{eq:defk})} below.
\end{prop}

\noindent{\sc Proof.} Let $n \in \{0,\dots,N-1\}$. Equation
(\ref{eq:1.10}) for $t=t_{n+1}$ reads
\begin{equation*}
\partial_t c(t_{n+1}) - D_c  \, \Delta c(t_{n+1})
+\bt\big(\U(t_{n+1}),c(t_{n+1})\big)
  =s_c-(s+\lambda)\, c(t_{n+1}).
\end{equation*}
We introduce the time discretization by setting 
\begin{equation*}
\label{eq:eqerr150}   R^{n+1}=  \Bigl( \frac{c(t_{n+1})-c(t_{n})}{k}
- c_t(t_{n+1}) \Bigr)
  +\bt\big(\U(t_{n})-\U(t_{n+1}),c(t_{n+1})\big).
\end{equation*}
We get
\begin{equation*}
 \frac{c(t_{n+1}) -c(t_n)}{k} - D_c \, \Delta c(t_{n+1})
+\bt\big(\U(t_n),c(t_{n+1})\big)
 =s_c-(s+\lambda) \, c(t_{n+1})+ R^{n+1}.
\end{equation*}
We apply $\Pi_{P_0}$ to this equation. By subtracting the result
from (\ref{eq:eqcd}) we get \stepcounter{err}
\begin{eqnarray}
\hspace{-1cm}&&  \Pi_{P_0} \Bigl( \frac{c(t_{n+1})-c(t_n)}{k} \Bigr)
 -\frac{c^{n+1}_h-c_h^n}{k}
  -D_c  \big( \Pi_{P_0} \Delta c(t_{n+1}) - \Delta_{h} c^{n+1}_h \big)
    \nonumber \\
\hspace{-1cm}&& +   \Pi_{P_0}  \bt\big(\U(t_n),c(t_{n+1})\big) -
\bt_h(\U^n_h,c^{n+1}_h) + \Pi_{P_0} \big( (s+\lambda) \,
c(t_{n+1})\big) - (s_h+\lambda) \, c^{n+1}_h  =\Pi_{P_0}R^{n+1}.
\label{eq:eqerr4}
\end{eqnarray}
We now introduce the discrete errors as follows.
Since $ c(t_{n+1}) - c(t_{n})=\int^{t_{n+1}}_{t_n} c_t(s)  \, ds$
one
checks that 
\begin{equation*}
 \Pi_{P_0} \Bigl( \frac{c(t_{n+1})-c(t_n)}{k} \Bigr)
 - \frac{c^{n+1}_h-c_h^n}{k}
 = \frac{1}{k} \int^{t_{n+1}}_{t_n}  \big(\Pi_{P_0} c_t(s) - \widetilde \Pi_{P_0} c_t(s)\big)  \, ds
 + \frac{1}{k} \, (\eps^{n+1}_{h,c}-\eps^{n}_{h,c}).
\end{equation*}
We also have 
\begin{equation*}
  \Pi_{P_0} \Delta c(t_{n+1}) - \Delta_{h} c^{n+1}_h
= \Pi_{P_0} \Delta c(t_{n+1}) -\Delta_{h} \big(\widetilde \Pi_{P_0}
c(t_{n+1})\big) +\Delta_{h} \eps^{n+1}_{h,c}. \label{eq:eqerr7}
\end{equation*}
Using the linearity of  $\bt_h$, one easily checks that
\begin{align*}
 \Pi_{P_0}  \bt\big(\U(t_n),c(t_{n+1})\big) -
\bt_h(\U^n_h,c^{n+1}_h) &= \bt_h(\U^n_h,\eps^{n+1}_{h,c})
 + \bt_h\big(\Eps{n}{h,\U},\widetilde
\Pi_{P_0} c(t_{n+1})\big)
 \nonumber \\
 &+ \Pi_{P_0} \bt\big(\U(t_n),c(t_{n+1})\big) -
\bt_h\big(\Pi_{\mathbf{RT_0}} \U(t_n),\widetilde \Pi_{P_0}
c(t_{n+1})\big)
 .
\label{eq:eqerr8}
\end{align*}
Lastly
 \stepcounter{err}
\begin{equation*}
  \Pi_{P_0} \big( (s+\lambda) \, c(t_{n+1}) \big) -
(s_h+\lambda) \, c^{n+1}_h  =    \Pi_{P_0} \big( (s+\lambda) \,
\eta^{n+1}_{h,c} \big)
 + (s_h+\lambda) \, \eps^{n+1}_{h,c}.
\end{equation*}
Using these relations  in (\ref{eq:eqerr4})
we get  (\ref{eq:eqerr26}). For any $m\in\{1,\dots,N\}$, the
consistency errors $C^m_{h,1} \in P_0$ and $C^m_{h,2} \in P_0$ are
given by 
\begin{eqnarray}
C^{m}_{h,1}&=&
 \Pi_{P_0}\Big( \frac{c(t_{m})-c(t_{m-1})}{k} -  c_t(t_{m})
+\bt\big(\U(t_{m-1})-\U(t_{m}),c(t_{m})\big) \Big)
\nonumber \\
&+& \Pi_{P_0} \big( (s+\lambda) \, \eta^{m}_{h,c} \big)- \frac{1}{k}
\int^{t_{m}}_{t_{m-1}} \big( \widetilde \Pi_{P_0}c_t(s)- \Pi_{P_0}
c_t(s) \big) \, ds,
\label{eq:deferrconsc1} \\
 C^{m}_{h,2}&=& D_c \, \Big( \Pi_{P_0} \Delta
c(t_{m}) - \Delta_{h} \big( \widetilde \Pi_{P_0} c(t_{m}) \big)
\Big) - \Big( \Pi_{P_0}\bt\big(\U(t_{m-1}),c(t_{m})\big) -
\bt_h(\Pi_{\mathbf{RT_0}} \U(t_{m-1}),\widetilde \Pi_{P_0} c(t_{m})
\big)\Big). \nonumber \label{eq:deferrconsc2}
\end{eqnarray}
A similar proof leads to (\ref{eq:eqerr27}) where the consistence
errors $\Theta^{m}_{h,1}\in P_0$ and $\Theta^{m}_{h,2}\in P_0$ are
defined
 for any $m\in\{1,\dots,N\}$ by
\stepcounter{err}
\begin{eqnarray}
\Theta^{m}_{h,1}&=&
  \Pi_{P_0}\Big( \frac{\theta(t_{m})-\theta(t_{m-1})}{k} - \theta_t(t_{m})
+\bt\big(\U(t_{m-1})-\U(t_{m}),\theta(t_{m})\big) \Big)
\nonumber \\
&+& \Pi_{P_0} ( s \, \eta^{m}_{h,\theta} )- \frac{1}{k}
\int^{t_{m}}_{t_{m-1}}   \big(\widetilde \Pi_{P_0}\theta_t(s)-
\Pi_{P_0} \theta_t(s)\big) \, ds,
\label{eq:deferrconstheta1} \\
\Theta^{m}_{h,2}&=& D_\theta \Big( \Pi_{P_0} \Delta \theta(t_{m}) -
\Delta_{h} \big( \widetilde \Pi_{P_0} \theta(t_{m}) \big) \Big) -
\Big( \Pi_{P_0}\bt\big(\U(t_{m-1}),\theta(t_{m})\big) -
\bt_h\big(\Pi_{\mathbf{RT_0}} \U(t_{m-1}),\widetilde \Pi_{P_0}
\theta(t_{m})\big) \Big). \nonumber \label{eq:deferrconstheta2}
\end{eqnarray}
We now consider the problem associated with the pressure. Let
$n\in\{0,\dots,N-1\}$ and $\psi_h \in P^{nc}_1$. Multiplying
equation (\ref{eq:1.12}) written for $t=t_{n+1}$ by $\psi_h$ and
integrating by parts, we get \stepcounter{err}
\begin{equation}
\label{eq:eqerr15}
  \big( \kappa(c(t_{n+1}),\theta(t_{n+1})) \nabla p(t_{n+1}),\nabla_h \psi_h \big)
= (\F(t_{n+1}),\nabla_h \psi_h)+(s,\psi_h).
\end{equation}
On the other hand, using (\ref{eq:eqpd}) and proposition
\ref{prop:propadjh}, we have
\begin{equation*}
  \big( \kappa(c^{n+1}_h,\theta^{n+1}_h)  \, \nabla_h p^{n+1}_h,\nabla_h \psi_h \big)
  =(\F^{n+1}_h,\nabla_h \psi_h)+(\widetilde \Pi_{P^{nc}_1} s_h,\psi_h).
\end{equation*}
Since $\nabla_h \psi_h \in P_0$, one checks that  $
(\F^{n+1}_h,\nabla_h \psi_h)=(\Pi_{\P_0} \F(t_{n+1}),\nabla_h
\psi_h)=(\F(t_{n+1}), \nabla_h\psi_h)$. According to
(\ref{eq:defprojp1nc}) we also have
 $(\widetilde \Pi_{P^{nc}_1}
s_h,\psi_h)=(s_h,\psi_h)$. Thus
\begin{equation*}
  \big( \kappa(c^{n+1}_h,\theta^{n+1}_h) \ \nabla_h p^{n+1}_h,\nabla_h \psi_h \big)
  =(\F(t_{n+1}),\nabla_h \psi_h)-(s_h,\psi_h).
\end{equation*}
Substracting (\ref{eq:eqerr15}) from the latter relation, we obtain
\stepcounter{err}
\begin{equation}
\label{eq:eqerr165} \big(\kappa(c(t_{n+1}),\theta(t_{n+1})) \,
\nabla p(t_{n+1})-\kappa(c^{n+1}_h,\theta^{n+1}_h) \,  \nabla_h
p^{n+1}_h,\nabla_h \psi_h \big) =-(s-s_h,\psi_h).
\end{equation}
We split the left-hand side as  
\begin{equation*}
 \kappa(c^{n+1}_h,\theta^{n+1}_h)  (\nabla p(t_{n+1})-\nabla_h p^{n+1}_h)
 +\big(\kappa(c(t_{n+1}),\theta(t_{n+1})) -
\kappa(c^{n+1}_h,\theta^{n+1}_h) \big)  \nabla p(t_{n+1}) .
\end{equation*}
Using a Taylor expansion, one can check that
\begin{equation*}
 \kappa(c(t_{n+1}),\theta(t_{n+1})) - \kappa(c^{n+1}_h,\theta^{n+1}_h)
 = (\eps^{n+1}_{h,c}+\eta^{n+1}_{h,c}) \,\kappa^{n+1}_{h,1} + (\eps^{n+1}_{h,\theta}+\eta^{n+1}_{h,\theta})\,
 \kappa^{n+1}_{h,2}.
\end{equation*}
We have set  for any $m\in\{0,\dots,N\}$ and $s\in [0,1]$
\stepcounter{err}
\begin{equation}
\left.
\begin{array}{lcl}
&&{\displaystyle
  c^{m}_h(s)=c^{m}_h+(c(t_{m})-c^{m}_h) \, s  , \hspace{1cm}
  \theta^{m}_h(s)=\theta^{m}_h+(\theta(t_{m})-\theta^{m}_h) \, s ,
  }\\
&&   {\displaystyle
  \kappa^{m}_{h,1}=\int^1_0 \kappa_x(c^{m}_h(s),\theta^{m}_h(s)) \, ds  ,
\hspace{1cm}
  \kappa^{m}_{h,2}=\int^1_0 \kappa_y(c^{m}_h(s),\theta^{m}_h(s)) \, ds.
  }
\end{array}
\right. \label{eq:defk}
\end{equation}
We also have
\begin{equation*}
  \nabla p(t_{n+1}) - \nabla_h p^{n+1}_h= \nabla_h \eps^{n+1}_{h,p}+\widetilde \nabla_h \eta^{n+1}_{h,p}.
\end{equation*}
Plugging these relations into (\ref{eq:eqerr165}) we get
(\ref{eq:eqerr28}). For all $m\in\{0,\dots,N\}$ we have
\stepcounter{err}
\begin{equation}
\label{eq:deferrconsp}
  \P^m_{h}=( \kappa^{m}_{h,1}   \eta^{m}_{h,c}+
  \kappa^{m}_{h,2} \eta^{m}_{h,\theta})\,  \nabla p(t_{m})
+\kappa(c^{m}_h,\theta^{m}_h)  \, \nabla_h \eta^{m}_{h,p}.
\end{equation}
We end with  the equation associated with $\U$. Let
$n\in\{0,\dots,N-1\}$. Applying the operator $\widetilde
\Pi_{\mathbf {RT_0}}$ to  (\ref{eq:1.12}) for $t=t_{n+1}$ we obtain
\begin{equation*}
\widetilde \Pi_{\mathbf {RT_0}} \U(t_{n+1})= \widetilde \Pi_{\mathbf
{RT_0}} \F(t_{n+1})- \widetilde \Pi_{\mathbf {RT_0}} \big( \kappa
(c(t_{n+1}),\theta(t_{n+1}) )
  \nabla p(t_{n+1})\big).
\end{equation*}
Let us substract this equation from (\ref{eq:equd}). Since
$\F_h^{n+1}=\Pi_{\P_0} \F(t_{n+1})$ we get 
\begin{eqnarray*}
\label{eq:eqerr23} \hspace{-.7cm}\widetilde \Pi_{\mathbf {RT_0}}
\U(t_{n+1})-\U^{n+1}_h
&=& \widetilde \Pi_{\mathbf {RT_0}} \big( \F(t_{n+1}) - \Pi_{\P_0} \F(t_{n+1})\big) \nonumber \\
&-&\widetilde \Pi_{\mathbf {RT_0}} \Big(
\kappa\big(c(t_{n+1}),\theta(t_{n+1}) \big)
 \, \nabla p(t_{n+1})- \kappa^{n+1}_h  \, \nabla_h p^{n+1}_h \Big).
\end{eqnarray*}
One easily checks that
\begin{equation*}
 \widetilde \Pi_{\mathbf {RT_0}} \U(t_{n+1})-\U^{n+1}_h  =\widetilde \Pi_{\mathbf {RT_0}} (
\U(t_{n+1}) - \Pi_{\mathbf{RT_0}} \U(t_{n+1}) ) +\Eps{n+1}{h,\U}.
\end{equation*}
\stepcounter{err}
Thus we get (\ref{eq:eqerr29}). For all $m\in \{0,\dots,N\}$, we
have
\begin{equation}
\label{eq:deferrconsu}
   \mathbf{U}^m_{h}= \widetilde \Pi_{\mathbf{RT_0}} \big(
 (\F(t_{m})  - \Pi_{\P_0}\F(t_{m})) -\Eta{m}{h,\U}-\P^m_{h} \big).
\end{equation}
This ends the proof of proposition \ref{propeqerr}. \qed


\subsection{Error estimates}

We first estimate the consistency errors.

\begin{prop}
\label{propesterrcons} For all $m\in\{1,\dots,N\}$ the consistency
errors satisfy\stepcounter{err}
\begin{eqnarray}
&& k \sum_{n=1}^m |C^n_{h,1}|^2 +  k \sum_{n=1}^m |\Theta^n_{h,1}|^2
 \le C \, (h^2+k^2)  ,
 \label{eq:esterrconsc}\\
\stepcounter{err} &&  k \sum_{n=1}^m \|C^n_{h,2}\|^2_{-1,h} + k
\sum_{n=1}^m \|\Theta^n_{h,2}\|^2_{-1,h}
 \le C \, h^2 \, ,
 \label{eq:esterrconsa}\\
\stepcounter{err} && |\P^m_h|+|\mathbf{U}^m_h| \le C \, h.
 \label{eq:esterrconsb}
\end{eqnarray}
\end{prop}
\noindent {\sc Proof.} Let $n\in\{1,\dots,N\}$. Since the operator
$\Pi_{P_0}$ is stable for  the $L^2$-norm we have
\begin{equation*}
 |\Pi_{P_0} R^n| \le |R^n| \le  \Bigl\vert
\frac{c(t_{n})-c(t_{n-1})}{k} - c_t(t_n) \Bigr\vert +\left|
\bt\big(\U(t_{n-1})-\U(t_n),c(t_n)\big)\right|.
\end{equation*}
Using a  Taylor expansion and the Cauchy-Schwarz inequality, we get
\begin{equation*}
 \Bigl\vert \frac{c(t_{n})-c(t_{n-1})}{k} - c_t(t_n) \Bigr\vert
  \le \frac{1}{k}\int^{t_{n-1}}_{t_n} |t_{n-1}-s| \, |c_{tt}(s)| \, ds
 \le \sqrt{k} \, \Bigl(  \int^{t_{n-1}}_{t_n}  |c_{tt}(s)|^2 \, ds \Bigr)^{1/2}.
\end{equation*}
On the other hand, since $\nabla c(t_n)|_{\partial \Omega}=0$, we
deduce from (\ref{eq:defbth})  by  integrating by parts
\begin{eqnarray*}
&&  \bt(\U(t_{n-1})-\U(t_n),c(t_n))= (\U(t_{n-1})-\U(t_n))\cdot
\nabla c(t_n).
\end{eqnarray*}
Using a Taylor expansion and  the Cauchy-Schwarz inequality, we get
\begin{eqnarray*}
&& \left|\left(\U(t_{n-1})-\U(t_n)\right)\cdot \nabla c(t_n)
\right|
  \le
 \sqrt{k} \, \|c\|_{L^\infty(0,T;H^1)} \left( \int^{t_{m-1}}_{t_m} \vert \U_t(s)\vert^2 \, ds \right)^{1/2}.
\end{eqnarray*}
Thus 
\begin{equation*}
\label{eq:esterrcons5}  |\Pi_{P_0} R^n|  \le  \sqrt{k} \,
  \Bigl(  \int^{t_{m-1}}_{t_m}  |c_{tt}(s)|^2 \, ds \Bigr)^{1/2}
 + \sqrt{k} \, \|c\|_{L^\infty(0,T;H^1)} \, \Bigl( \int^{t_{m-1}}_{t_m} \vert \U_t(s)\vert^2 \, ds
 \Bigr)^{1/2}.
\end{equation*}
Thanks to the stability of $\Pi_{P_0}$ for the $L^2$-norm
and to (\ref{eq:esterrintc}) we have 
\begin{equation*}
  \left| \Pi_{P_0} \big( (s+\lambda) \, \eta^n_{h,c} \big)\right|
  \le  C \, h \, (\Vert s \Vert_{L^\infty(0,T;L^2)}+\lambda)\, \|c\|_{L^\infty(0,T;H^1)}.
\label{eq:esterrcons9}
\end{equation*}
The Cauchy-Schwarz inequality and (\ref{eq:esti}) allow
us to write 
\begin{equation*}
  \int^{t_m}_{t_{m-1}} \bigl\vert \Pi_{P_0}  c_t(s) - \widetilde \Pi_{P_0} c_t(s) \bigr\vert \, ds
\le C  \, h  \, \sqrt{k} \,  \Bigl( \int^{t_m}_{t_{m-1}}
\|c_t(s)\|^2_{1+r} \, ds \Bigr)^{1/2}. \label{eq:esterrcons8}
\end{equation*}
By plugging these  estimates into  definition
(\ref{eq:deferrconsc1}) we get
 \begin{eqnarray*}
 k \,  |C^n_{h,1}|^2 &\le& k^2 \, \|c\|^2_{L^\infty(0,T;H^1)}
  \int^{t_n}_{t_{n-1}} \vert \U_t(s)\vert^2 \, ds +k^2 \int^{t_n}_{t_{n-1}} |c_{tt}(s)|^2 \, ds
\\
&+&C \, h^2 \int^{t_n}_{t_{n-1}} \|c_t(s)\|^2_{1+r} \, ds+  C \, k
\, h^2 \, \|c\|^2_{L^\infty(0,T;H^1)}.
 \end{eqnarray*}
Summing up the latter relation for  $n=1$ to $m\in\{1,\dots,N\}$ and
using  a similar work on $\Theta^m_{h,1}$ we get
(\ref{eq:esterrconsc}). Now let $n\in\{1,\dots,N\}$. Using
propositions \ref{propconslap} and \ref{propconsbh} we have
\begin{equation*}
\label{eq:esterrcons6}
  \|\Pi_{P_0} \Delta c(t_n) -\Delta_{h}\big(\widetilde \Pi_{P_0} c(t_n)\big)\|_{-1,h}
\le C \, h\, \|c\|_{L^\infty(0,T;H^2)}
\end{equation*}
and
\begin{equation*}
\|\Pi_{P_0}\bt(\U(t_{n-1}),c(t_n)) - \bt_h(\Pi_{\mathbf{RT_0}}
\U(t_{n-1}), \widetilde \Pi_{P_0} c(t_n))\|_{-1,h}   \le C\,  h \,
\|c\|_{L^\infty(0,T;H^1)} \|\U\|_{L^\infty(0,T;\H^{1+s})}.
\end{equation*}
Plugging these estimates into  definition (\ref{eq:deferrconsc2})
and summing up from $n=1$ to $m$, we obtain
\begin{equation*}
   k \sum_{n=1}^m \|C^n_{h,2}\|^2_{-1,h} \le
  C \, h^2 \, \bigl( \|c\|^2_{L^\infty(0,T;H^1)} \, \|\U\|^2_{L^\infty(0,T;\H^{1+s})}
  +  \|c\|^2_{L^\infty(0,T;H^2)} \bigr).
\end{equation*}
A similar work on $\Theta^m_{h,2}$ then leads to
(\ref{eq:esterrconsa}). We finally prove (\ref{eq:esterrconsb}). Let
$m\in\{1,\dots,N\}$. On the one hand, we have by
(\ref{eq:deferrconsp})
 \begin{equation*}
  |\P^m_h|
  \le
  ( \|\kappa^m_{h,1}\|_{\infty}  \, |\eta^m_{h,c}|
  +\|\kappa^m_{h,2}\|_{\infty}   \, |\eta^m_{h,\theta}| ) \, |\nabla
  p(t_m)|
  +\|\kappa(c^m_h,\theta^m_h)\|_{\infty}  |\nabla_h \eta^m_{h,p}|.
 \end{equation*}
Using estimates (\ref{eq:esterrintc})--(\ref{eq:esterrintu}) we get 
\begin{equation*}
\label{eq:esterrcons19}
  |\P^m_h| \le C \, h \, \|\kappa\|_{W^{1,\infty}((0,1) \times (0,\infty))}   (\|c\|_{L^\infty(0,T;H^1)}
  +\|\theta\|_{L^\infty(0,T;H^1)}
  + \|p\|_{L^\infty(0,T;H^2)}).
\end{equation*}
On the other hand definition (\ref{eq:deferrconsu})  leads to
\begin{equation*}
 |\mathbf{U}^m_h| \le  \bigl\vert  \widetilde \Pi_{\mathbf{RT_0}} \big( \F(t_m)-\Pi_{\P_0} \F(t_m)\big) \bigr\vert
  +|\Eta{m}{h,\U}|+|\P^m_h|.
\end{equation*}
Using the stability of $\widetilde \Pi_{\mathbf{RT_0}}$ for the
$\L^2$-norm and proposition \ref{prop:errint} we have
\begin{equation*}
  \left| \widetilde \Pi_{\mathbf{RT_0}} \big( \F(t_m)-\Pi_{\P_0} \F(t_m)\big) \right| \le
  |\F(t_m)-\Pi_{\P_0} \F(t_m)|
   \le C \, h \, \|\F\|_{L^\infty(0,T;\H^1)}.
\end{equation*}
Using moreover  (\ref{eq:esterrintu})  we obtain
\begin{align*}
\label{eq:esterrcons20} |\mathbf{U}^m_h|
  &\le  C  \, h \,  \|\kappa\|_{W^{1,\infty}((0,1) \times (0,\infty))} (\|c\|_{L^\infty(0,T;H^1)}
  +\|\theta\|_{L^\infty(0,T;H^1)}) \\ &+ C \, (\|p\|_{L^\infty(0,T;H^2)}
  +  \|\F\|_{L^\infty(0,T;\H^1)}+\|\U\|_{L^\infty(0,T;\H^1)}).
\end{align*}
We have proven (\ref{eq:esterrconsb}). \qed


\vspace{.2cm}

\noindent Using the former proposition  we are now able to estimate
the discrete errors.
\begin{prop}
\label{prop:esterrd} There exists some real  $k_0>0$ such that for
any  $k<k_0$ and $m\in\{1,\dots,N\}$ \stepcounter{err}
  \begin{eqnarray}
&&  |\eps^m_{h,c}|^2+|\eps^m_{h,\theta}|^2+k \sum_{n=1}^m
\left(\|\eps^m_{h,c}\|^2_h
    +\|\eps^m_{h,\theta}\|^2_h\right) \le C \, (h^2+k^2),
\label{eq:esterrda} \\
\stepcounter{err} &&     |\nabla_h \eps^m_{h,p}|+ |\eps^m_{h,\U}|
\le C \, (h+k).
 \label{eq:esterrdb}
  \end{eqnarray}
\end{prop}
\noindent {\sc Proof.} Multiplying (\ref{eq:eqerr26}) by $2  \, k \,
\eps^{n+1}_{h,c}$, we obtain \stepcounter{err}
\begin{eqnarray}
\label{eq:esterrd1} &&  \hspace{-1.5cm}\Bigl(
\frac{\eps^{n+1}_{h,c}-\eps^n_{h,c}}{k}, \,  2 \,  k \,
\eps^{n+1}_{h,c}\Bigr)
  - 2 \,  D_c \, k\, (\Delta_{h} \eps^{n+1}_{h,c},\eps^{n+1}_{h,c})
  + 2  \, k \, \b_h\big(\U^n_h,\eps^{n+1}_{h,c},\eps^{n+1}_{h,c}\big)+2  \, k  \, \lambda \, |\eps^{n+1}_{h,c}|^2
  \nonumber \\
&& \hspace{1.5cm}
 =2  k \, (C^{n+1}_{h+1}+C^{n+1}_{h,2},\eps^{n+1}_{h,c})
 -2  k\, ( s_h,\vert \eps^{n+1}_{h,c} \vert^2)
 -2  k \, \b_h\big(\Eps{n}{h,\U},\widetilde \Pi_{P_0} c(t_{n+1}),\eps^{n+1}_{h,c}\big).
\end{eqnarray}
Using an algebraic identity we have 
\begin{equation*}
\label{eq:esterrd2}  \Bigl( \frac{\eps^{n+1}_{h,c}-\eps^n_{h,c}}{k},
2 \, k \, \eps^{n+1}_{h,c}\Bigr) =
|\eps^{n+1}_{h,c}|^2-|\eps^n_{h,c}|^2+|\eps^{n+1}_{h,c}-\eps^n_{h,c}|^2.
\end{equation*}
We know by propositions \ref{prop:posbh} and \ref{prop:coerlap} that
\begin{equation*}
\label{eq:esterrd3}
  - 2 \,  k\, (\Delta_{h} \eps^{n+1}_{h,c},\eps^{n+1}_{h,c})
  = 2  \, k \,  \|\eps^{n+1}_{h,c}\|^2_h \, , \hspace{2cm}
  \b_h\big(\U^n_h,\eps^{n+1}_{h,c},\eps^{n+1}_{h,c}\big) \ge 0.
\end{equation*}
We have
\begin{equation*}
  2 \, k \, (s_h,|\eps^{n+1}_{h,c}|^2) \le 2 \, k \,
  \|s_h\|_{\infty} \, |\eps^{n+1}_{h,c}|^2 \le C \, k\,
  |\eps^{n+1}_{h,c}|^2.
\end{equation*}
 Using the  Young inequality, we also write 
\begin{eqnarray*}
&& \left| 2  k \, (C^{n+1}_{h,1},\eps^{n+1}_{h,c})\right| \le 2  k\,
|C^{n+1}_{h,1}| \, |\eps^{n+1}_{h,c}| \le k \, \lambda \,
|\eps^{n+1}_{h,c}|^2+   C \, k \, |C^{n+1}_{h,1}|^2,
\label{eq:esterrd65} \\ 
&& |2  k \, (C^{n+1}_{h,2},\eps^{n+1}_{h,c})| \le 2  k\,
\|C^{n+1}_{h,2}\|_{-1,h} \, \|\eps^{n+1}_{h,c}\|_h \le D_c \,
\frac{k}{2} \,  \|\eps^{n+1}_{h,c}\|^2_h+ C \, k  \,
\|C^{n+1}_{h,2}\|^2_{-1,h}. \label{eq:esterrd6}
\end{eqnarray*}
We are left with the term $\b_h\big(\Eps{n}{h,\U},\widetilde
\Pi_{P_0} c(t_{n+1}),\eps^{n+1}_{h,c}\big)$. We have  $\Eps{n}{h,\U}
\in \mathbf{RT_0}$. Using the divergence formula, proposition
\ref{prop:propdivsh} and (\ref{eq:defprojrt0}), one easily checks
that $\hbox{div} \, \Eps{n}{h,\U}=0$. Thus we can apply proposition
\ref{prop:stabbth} to get \stepcounter{err}
\begin{equation}
\label{eq:esterrd79} \big|\b_h(\Eps{n}{h,\U},\widetilde \Pi_{P_0}
c(t_{n+1}),\eps^{n+1}_{h,c})\big| \le C \, |\Eps{n}{h,\U}| \,
\|\widetilde \Pi_{P_0} c(t_{n+1})\|_h \,
 \|\eps^{n+1}_{h,c}\|_h.
\end{equation}
Let us first bound $\|\widetilde \Pi_{P_0} c(t_{n+1})\|_h$. We have
\begin{equation*}
\|\widetilde \Pi_{P_0} c(t_{n+1})\|_h \le \|\widetilde \Pi_{P_0}
c(t_{n+1})-\Pi_{P_0} c(t_{n+1})\|_h+\|\Pi_{P_0} c(t_{n+1})\|_h.
\end{equation*}
Using an inverse inequality (see proposition 1.2 in \cite{zimconv})
and (\ref{eq:esti})
\begin{equation*}
\|\widetilde \Pi_{P_0} c(t_{n+1})-\Pi_{P_0} c(t_{n+1})\|_h \le
\frac{C}{h} \, |\widetilde \Pi_{P_0} c(t_{n+1})-\Pi_{P_0}
c(t_{n+1})| \le C.
\end{equation*}
Moreover, according to \cite{eym} (p. 776), we have $\|\Pi_{P_0}
c(t_{n+1})\|_h \le C \, \|c\|_{L^\infty(0,T;H^1)}$. Thus
$\|\widetilde \Pi_{P_0} c(t_{n+1})\|_h  \le C$.
We now estimate $|\Eps{n}{h,\U}|$. Using the stability of
$\widetilde \Pi_{\mathbf{RT_0}}$ for the $\L^2$-norm
 and the  Cauchy-Schwarz inequality in  (\ref{eq:eqerr29}) we get
 \stepcounter{err}
\begin{equation}
\label{eq:esterrd95}
  |\Eps{n}{h,\U}| \le
  \|\kappa\|_{W^{1,\infty}((0,1) \times (0,\infty))} \, \Big( \|\nabla p\|_{L^\infty(0,T;\L^2)} \,
  (|\eps^{n+1}_{h,c}|+|\eps^{n+1}_{h,\theta}|)
  +|\nabla_h \eps^{n+1}_{h,p}| \Big)
  +|\mathbf{U}^{n+1}_h|.
\end{equation}
We  bound $\eps^{n+1}_{h,p}$ as follows. Setting
$\psi_h=\eps^{n+1}_{h,p}$ in (\ref{eq:eqerr28}) and using  the
Cauchy-Schwarz inequality, we get
\begin{eqnarray*}
\label{eq:esterrd10}
 \big( \kappa(c^{n+1}_h,\theta^{n+1}_h)  \, \nabla_h \eps^{n+1}_{h,p}, \nabla_h \eps^{n+1}_{h,p}\big)
 &\le& \|\kappa\|_{W^{1,\infty}((0,1) \times (0,\infty))} \, \|\nabla
 p\|_{L^\infty(0,T;\L^2)} \,
  ( |\eps^{n+1}_{h,c}|+|\eps^{n+1}_{h,\theta}| )  \, |\nabla_h \eps^{n+1}_{h,p}|
  \nonumber \\
   &+& |\P^{n+1}_h|  \, |\nabla_h \eps^{n+1}_{h,p}|
   +|s-s_h| \, |\eps^{n+1}_{h,p}|.
\end{eqnarray*}
The left-hand side  is such that $ (
\kappa(c^{n+1}_h,\theta^{n+1}_h)  \nabla_h \eps^{n+1}_{h,p},
\nabla_h \eps^{n+1}_{h,p} ) \ge \kappa_{inf} \,  | \nabla_h
\eps^{n+1}_{h,p}|^2$. As for the right-hand side, we have
$s_h=\Pi_{P_0} s$ and $\eps^{n+1}_{h,p} \in P^{nc}_1 \cap L^2_0$.
Thus, according to propositions  \ref{proppoinp1nc} and
\ref{prop:errint}
\begin{equation*}
 |s-s_h| \, |\eps^{n+1}_{h,p}| \le C \, h \, \|s\|_{L^\infty(0,T;H^1)} \, |\nabla_h
 \eps^{n+1}_{h,p}|.
\end{equation*}
Finally  $  |\P^{n+1}_h| \le C \, h$ thanks to
(\ref{eq:esterrconsb}). Therefore we obtain \stepcounter{err}
\begin{equation}
\label{eq:esterrd105}
  |\nabla_h \eps^{n+1}_{h,p}| \le C \, (h+|\eps^{n+1}_{h,c}|+|\eps^{n+1}_{h,\theta}| ).
\end{equation}
Let us plug this estimate into (\ref{eq:esterrd95}). Since $
|\mathbf{U}^{n+1}_h| \le C \, h$ thanks to (\ref{eq:esterrconsb}),
we get
 \stepcounter{err}
\begin{equation}
\label{eq:esterrd11}
  |\Eps{n}{h,\U}| \le C \, (h+|\eps^{n+1}_{h,c}|+|\eps^{n+1}_{h,\theta}| ).
\end{equation}
Now, plugging this bound into (\ref{eq:esterrd95}) and using the
Young inequality, we obtain
\begin{eqnarray*}
 k\, \big|\b_h(\Eps{n}{h,\U},\widetilde \Pi_{P_0} c(t_{n+1}),\eps^{n+1}_{h,c})\big|
& \le &C \, k \, (h+|\eps^{n+1}_{h,c}|+|\eps^{n+1}_{h,\theta}| ) \,
\|\eps^{n+1}_{h,c}\|_h
\nonumber \\
&  \le&  D_c \, \frac{k}{2} \, \|\eps^{n+1}_{h,c}\|^2_h + C \, k\,
(h^2+|\eps^{n+1}_{h,c}|^2+|\eps^{n+1}_{h,\theta}|^2 ).
  \label{eq:esterrd12}
\end{eqnarray*}
Now we have treated all the terms in (\ref{eq:esterrd1}). This
equation implies
\begin{equation*}
|\eps^{n+1}_{h,c}|^2-|\eps^n_{h,c}|^2 +D_c \, k \,
\|\eps^{n+1}_h\|^2_h
  \le C \, k  \, \big(h^2+|\eps^{n+1}_{h,c}|^2+|\eps^{n+1}_{h,\theta}|^2 +|C^{n+1}_{h,1}|^2
  +\|C^{n+1}_{h,2}\|^2_{-1,h}\big).
\end{equation*}
Let $m\in \{1,\dots,N\}$. Let us sum up the latter estimate from
$n=0$ to $m-1$.
Thanks to (\ref{prop:errint})
\begin{equation*}
|\eps^0_{h,c}| =|\widetilde \Pi_{P_0} c_0-c^0_h|=|\widetilde
\Pi_{P_0} c_0-\Pi_{P_0} c_0|
 \le C  \, h \, \|c\|_{L^\infty(0,T;H^2)}.
\end{equation*}
Using moreover estimates (\ref{eq:esterrconsc}) and
(\ref{eq:esterrconsa}) we get
\begin{equation*}
|\eps^m_{h,c}|^2+ D_c \, k   \sum_{n=1}^m \|\eps^{n}_{h,c}\|^2_h \le
C \, k
 \sum_{n=1}^m (|\eps^{n}_{h,c}|^2+|\eps^{n}_{h,\theta}|^2 ) +C \,
(h^2+k^2).
\end{equation*}
Summing this relation with the one obtained by a similar work on
(\ref{eq:eqerr27}) we obtain
\begin{equation*}
 |\eps^m_{h,c}|^2+ |\eps^m_{h,\theta}|^2+k  \sum_{n=1}^m ( D_c \, \|\eps^{n}_{h,c}\|^2_h +D_\theta \, \|\eps^{n}_{h,\theta}\|^2_h)
   \le  C \, k    \sum_{n=1}^m (|\eps^{n}_{h,c}|^2+|\eps^{n}_{h,\theta}|^2 )
  + C   \, (h^2+k^2).
\end{equation*}
Using a discrete Gronwall lemma (see lemma 5.2 in \cite{zimconv}) we
get (\ref{eq:esterrda}).
 Then (\ref{eq:esterrd105}) and (\ref{eq:esterrd11}) imply (\ref{eq:esterrdb}). \qed

\vspace{.2cm}

\noindent By combining proposition \ref{prop:esterrd} with estimates
(\ref{eq:esterrintc})-(\ref{eq:esterrintu}),
 we obtain  finally the following result.

\begin{theor}
There exists a real $k_0>0$ such that for all $k<k_0$ and
$m\in\{1,\dots,N\}$ 
  \begin{eqnarray*}
    &&|e^m_{h,c}|^2+|e^m_{h,\theta}|^2+k \sum_{n=1}^m \big(\|\widetilde \Pi_{P_0} e^m_{h,c}\|^2_h
    +\|\widetilde \Pi_{P_0} e^m_{h,\theta}\|^2_h\big) \le C \,
    (h^2+k^2) \, , \\
&&    |\widetilde \nabla_h e^m_{h,p}|+ |e^m_{h,\U}| \le C  \, (h+k).
  \end{eqnarray*}
\end{theor}



\bibliographystyle{plain}



%

\end{document}